\documentclass[12pt]{article}

\usepackage{amsfonts,amssymb,amsmath,amsthm,epsfig,euscript}

\theoremstyle{definition}
{

\newtheorem{problem}{Problem}
}

\newtheorem{theorem}{Theorem}
\newtheorem{lemma}[theorem]{Lemma}

\long\def\symbolfootnote[#1]#2{\begingroup
\def\thefootnote{\fnsymbol{footnote}}\footnote[#1]{#2}\endgroup}

\newcommand{\coinv}[1][\sigma]{\mathrm{coinv}(#1)}
\newcommand{\coinvi}[1][\sigma^{(i)}]{\mathrm{coinv}(#1)}

\newcommand{\red}[1][\sigma]{\mathrm{red}(#1)}
\newcommand{\inv}[1][\sigma]{\mathrm{inv}(#1)}

\newcommand{\qbin}[3]{\genfrac{[}{]}{0pt}{}{#1}{#2}_{#3}}

\newcommand{\sg}{\sigma}

\newcommand{\RLmax}[1][\sigma]{\mathrm{RLmax}(#1)}

\newcommand{\fig}[2]{\begin{figure}[ht]
\centerline{\scalebox{.66}{\epsfig{file=#1.eps}}}
\caption{#2}
\label{fig:#1}
\end{figure}}

\title{Simple marked mesh patterns}

\author{
Sergey Kitaev \\
\small University of Strathclyde\\[-0.8ex]
\small Livingstone Tower, 26 Richmond Street\\[-0.8ex]
\small Glasgow G1 1XH, United Kingdom\\[-0.8ex]
\small \texttt{sergey.kitaev@cis.strath.ac.uk}
\and
Jeffrey Remmel \\
\small Department of Mathematics\\[-0.8ex]
\small University of California, San Diego\\[-0.8ex]
\small La Jolla, CA 92093-0112. USA\\[-0.8ex]
\small \texttt{remmel@math.ucsd.edu}
}

\date{\small Submitted: Date 1;  Accepted: Date 2;
 Published: Date 3.\\
\small MR Subject Classifications: 05A15, 05E05}

\begin{document}
\maketitle

\begin{abstract}
\noindent \

In this paper we begin the first systematic study of distributions of simple marked mesh patterns. Mesh patterns were introduced recently by Br\"and\'en and Claesson \cite{BrCl} in connection  with permutation statistics. We provide explicit generating functions in several general cases, and develop recursions to compute the numbers in question in some other cases. Certain $q$-analogues are discussed. Moreover, we consider two modifications of the notion of a marked mesh pattern and provide enumerative results for them. \\    

\end{abstract}

\tableofcontents

\section{Introduction}

The notion of mesh patterns was introduced by Br\"and\'en and Claesson \cite{BrCl} to provide explicit expansions for certain permutation statistics as, possibly infinite, linear combinations of (classical) permutation patterns.  This notion was further studied in \cite{HilJonSigVid,Ulf}. In particular, the notion of a mesh pattern was extended to that of a {\em marked mesh pattern} by \'Ulfarsson in \cite{Ulf}. 

In this paper, we study the number of occurrences of what 
we call {\em simple marked mesh patterns}. Let 
$\mathbb{N} = \{0,1,2, \ldots \}$ denote the set of 
natural numbers and $S_n$ denote the symmetric group of permutations 
 of $1, \ldots, n$. 
If $\sigma = \sg_1 \ldots \sg_n \in S_n$, then we will consider the 
graph of $\sg$, $G(\sg)$, to be the set of points $(i,\sg_i)$ for 
$i =1, \ldots, n$.  For example, the graph of the permutation 
$\sg = 471569283$ is pictured in Figure 
\ref{fig:basic}.  Then if we draw a coordinate system centered at a 
point $(i,\sg_i)$, we will be interested in  the points that 
lie in the four quadrants I, II, III, and IV of that 
coordinate system as pictured 
in Figure \ref{fig:basic}.  For any $a,b,c,d \in  
\mathbb{N}$  and any $\sg = \sg_1 \ldots \sg_n \in S_n$, 
we say that $\sg_i$ matches the 
simple marked mesh pattern $MMP(a,b,c,d)$ in $\sg$ if  in $G(\sg)$ relative 
to the coordinate system which has the point $(i,\sg_i)$ as its  
origin, there are 
$\geq a$ points in quadrant I, 
$\geq b$ points in quadrant II, $\geq c$ points in quadrant 
III, and $\geq d$ points in quadrant IV.  
For example, 
if $\sg = 471569283$, the point $\sg_4 =5$  matches 
the simple marked mesh pattern $MMP(2,1,2,1)$ since relative 
to the coordinate system with origin $(4,5)$,  
there are 3 points in $G(\sg)$ in quadrant I, there is 
1 point in $G(\sg)$ in quadrant II, there are 2 points in $G(\sg)$ in quadrant III, and there are 2 points in $G(\sg)$ in 
quadrant IV.  Note that if a coordinate 
in $MMP(a,b,c,d)$ is 0, then there is no condition imposed 
on the points in the corresponding quadrant. In addition, we shall 
consider patterns  $MMP(a,b,c,d)$ where 
$a,b,c,d \in \mathbb{N} \cup \{\emptyset\}$. Here when 
a coordinate of $MMP(a,b,c,d)$ is the empty set, then for $\sg_i$ to match  
$MMP(a,b,c,d)$ in $\sg = \sg_1 \ldots \sg_n \in S_n$, 
it must be the case that there are no points in $G(\sg)$ relative 
to the coordinate system with origin $(i,\sg_i)$ in the corresponding 
quadrant. For example, if $\sg = 471569283$, the point 
$\sg_3 =1$ matches 
the marked mesh pattern $MMP(4,2,\emptyset,\emptyset)$ since relative 
to the coordinate system with origin $(3,1)$, 
there are 6 points in $G(\sg)$ in quadrant I, 
2 points in $G(\sg)$ in quadrant II, no  points in $G(\sg)$ in quadrant III, and no  points in $G(\sg)$ in quadrant IV.  We let 
$mmp^{(a,b,c,d)}(\sg)$ denote the number of $i$ such that 
$\sg_i$ matches the marked mesh pattern $MMP(a,b,c,d)$ in $\sg$.

\fig{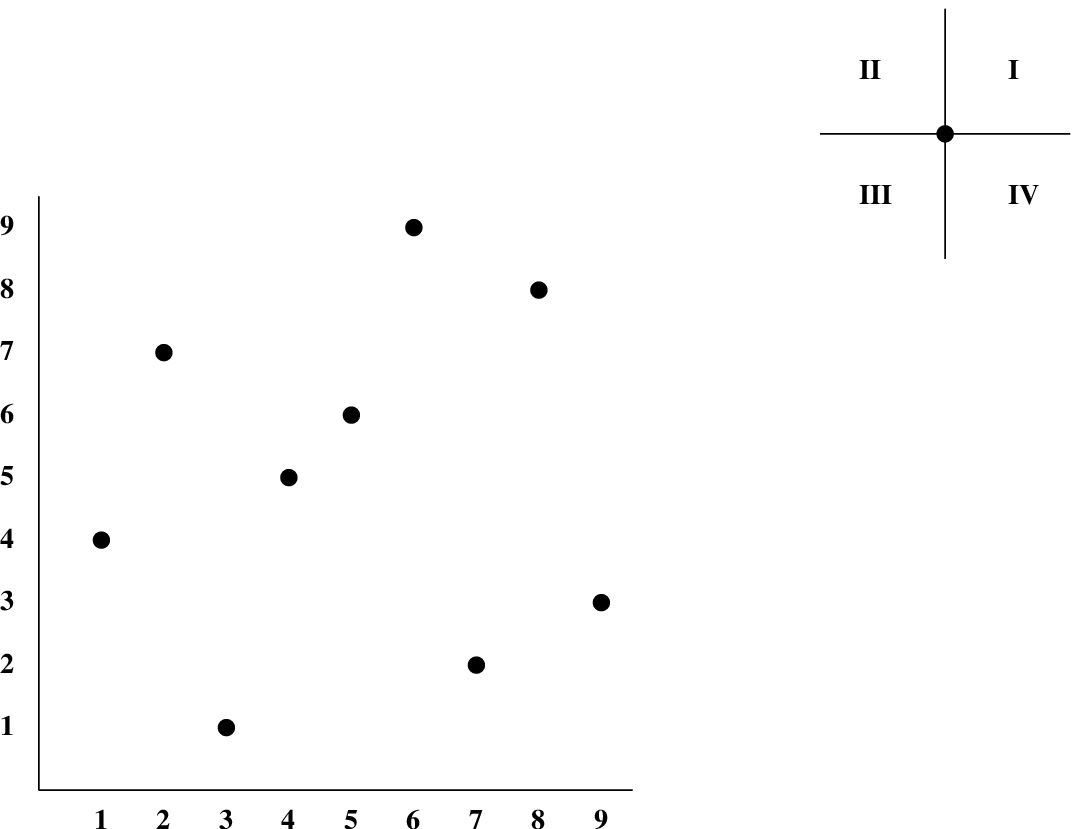}{The graph of $\sg = 471569283$.}

Given a sequence $\sg = \sg_1 \ldots \sg_n$ of distinct integers,
let $\red[\sg]$ be the permutation found by replacing the
$i$-th largest integer that appears in $\sg$ by $i$.  For
example, if $\sg = 2754$, then $\red[\sg] = 1432$.  Given a
permutation $\tau=\tau_1 \ldots \tau_j\in S_j$, we say that the pattern $\tau$ {\em occurs} in $\sg = \sg_1 \ldots \sg_n \in S_n$ provided   there exists 
$1 \leq i_1 < \cdots < i_j \leq n$ such that 
$\red[\sg_{i_1} \ldots \sg_{i_j}] = \tau$.   We say 
that a permutation $\sg$ {\em avoids} the pattern $\tau$ if $\tau$ does not 
occur in $\sg$. Let $S_n(\tau)$ denote the set of permutations in $S_n$ 
which avoid $\tau$. In the theory of permutation patterns,  $\tau$ is called a {\em classical pattern}. See \cite{kit} for a comprehensive introduction to 
permutation patterns.

Given a permutation $\sg = \sg_1 \ldots \sg_n \in S_n$, 
we say that $\sg_j$ is a \emph{right-to-left maximum} of $\sg$ if $\sg_j > \sg_i$ for 
all $i>j$. We let $\RLmax[\sg]$ denote the number of 
right-to-left maxima of~$\sg$. 

The main goal of this paper is to study  the generating 
functions 
\begin{equation} \label{Rabcd}
R^{(a,b,c,d)}(t,x) = 1 + \sum_{n\geq 1} \frac{t^n}{n!} 
\sum_{\sg \in S_n} x^{mmp^{(a,b,c,d)}(\sg)}.
\end{equation}

For any $a,b,c,d \in \{\emptyset\} \cup \mathbb{N}$, let 
\begin{equation} \label{Rabcdn}
R^{(a,b,c,d)}_n(x) = \sum_{\sg \in S_n} x^{mmp^{(a,b,c,d)}(\sg)}.
\end{equation}
Note that there is a natural action of the symmetries of the square on  
the graphs of permutations.  That is, one can identify each permutation 
$\sg \in S_n$ with a permutation matrix $P(\sg)$ and it is 
clear that rotating such matrices by 90 degrees counterclockwise or reflecting 
permutation matrices about any central axis or any diagonal axis 
preserves the property of being a 
permutation matrix. For any permutation $\sg = \sg_1 \ldots \sg_n \in S_n$, 
let 
$\sg^r = \sg_n \ldots \sg_1$ be the {\em reverse} of $\sg$ and 
$\sg^c = (n+1 -\sg_1) \ldots (n+1 - \sg_n)$ be the {\em complement} 
of $\sg$.  Then it is easy to see that the map 
$\sg \rightarrow \sg^r$ corresponds 
to reflecting a permutation matrix about its central vertical 
axis and the map $\sg \rightarrow \sg^c$ corresponds to 
reflecting a permutation matrix  about its central horizontal axis. 
It is also easy to see that 
rotating a permutation 
matrix by 90 degrees corresponds to replacing $\sg$ by $(\sg^{-1})^r$.
It follows that the map $\sg \rightarrow (\sg^{-1})^r$ shows that 
$R_n^{(a,b,c,d)}(x) = R_n^{(d,a,b,c)}(x)$, 
the map $\sg \rightarrow \sg^r$ shows 
that for $R_n^{(a,b,c,d)}(x) = R_n^{(b,a,d,c)}(x)$ and the 
map $\sg \rightarrow \sg^c$ shows 
that for $R_n^{(a,b,c,d)}(x) = R_n^{(d,c,b,a)}(x)$. Moreover, taking the inverse of permutations shows that $R_n^{(a,b,c,d)}(x) = R_n^{(a,d,c,b)}(x)$, applying the composition of reverse and complement shows that $R_n^{(a,b,c,d)}(x) = R_n^{(c,d,a,b)}(x)$, applying the composition of reverse and inverse shows that $R_n^{(a,b,c,d)}(x) = R_n^{(b,c,d,a)}(x)$, and finally, applying the composition of reverse, complement and inverse,  that is, the map $\sg \rightarrow (\sg^{rc})^{-1}$, shows that  $R_n^{(a,b,c,d)}(x) = R_n^{(c,b,a,d)}(x)$. Any other composition of the three trivial bijections on $S_n$ will not give us any new symmetries.  Thus we have the following lemma:
\begin{lemma} \label{sym}
For any $a,b,c,d \in \{\emptyset\} \cup \mathbb{N}$, 
\begin{eqnarray*}
R_n^{(a,b,c,d)}(x) &=& R_n^{(d,a,b,c)}(x) = R_n^{(c,b,a,d)}(x) =\\
R_n^{(b,a,d,c)}(x) &=& R_n^{(d,c,b,a)}(x) = R_n^{(a,d,c,b)}(x) = R_n^{(c,d,a,b)}(x) = R_n^{(b,c,d,a)}(x).   
\end{eqnarray*}
\end{lemma}

\section{$R_n^{(k,0,0,0)}(x)=R_n^{(0,k,0,0)}(x)=R_n^{(0,0,k,0)}(x)=R_n^{(0,0,0,k)}(x)$ and $R_n^{(=k,0,0,0)}(x)$}

The equalities in the section title are true by Lemma \ref{sym}. Thus we only need to consider the simple marked mesh pattern 
$MMP(k,0,0,0)$. 

First assume that $k \geq 1$. 
It is easy to see that $R^{(k,0,0,0)}_n(x) = n!$ if 
$n \leq k$.  For $n > k$, suppose that we start 
with a permutation $\sg  = \sg_1 \ldots \sg_n \in S_n$ and then 
we let $\sg^{(i)}$ denote the permutation of $S_{n+1}$ 
that results by  adding $1$ to each element 
of $\sg$ and inserting $1$ before the element $\sg_i+1$ if 
$1\leq i \leq n$ and inserting 1 at the end if $i =n+1$.  Then 
it is easy to see that 
\begin{equation}\label{eq:k1}
mmp^{(k,0,0,0)}(\sg^{(i)}) = \begin{cases} 
mmp^{(k,0,0,0)}(\sg) & \mbox{if $i > n-k$}, \\
1+ mmp^{(k,0,0,0)}(\sg) & \mbox{if $i \leq  n-k$}.
\end{cases}
\end{equation}
That is, clearly the placement of 1 in position $i$ in $\sg^{(i)}$ 
cannot effect the elements 
in quadrant I relative to any pair $(j,\sg_j+1)$ for 
any $j <i$  or relative to any pair $(j+1,\sg_j+1)$ for any $j \geq i$. 
In addition, the $1$ in position $i$ can contribute 
to $mmp^{(k,0,0,0)}(\sg^{(i)})$ if and only if $i \leq n+1-k$. 
It then follows that for $n \geq k$, 
\begin{equation}\label{eq:k2}
R^{(k,0,0,0)}_{n+1}(x) = (k + x(n+1-k))R^{(k,0,0,0)}_n(x).
\end{equation}
Iterating this recursion, we see that 
\begin{equation}\label{eq:k3}
R^{(k,0,0,0)}_{k+s}(x) = k! \prod_{i=1}^s(k+ix) 
\end{equation}
for all $s \geq 1$. In this case, we can form a simple generating 
function. That is, let 
\begin{equation} \label{Rabcd-}
P^{(k,0,0,0)}(t,x) = \sum_{n\geq k} \frac{t^{n-k}}{(n-k)!} 
\sum_{\sg \in S_n} x^{mmp^{(k,0,0,0)}(\sg)}.
\end{equation}
Then we know that 
\begin{eqnarray*}
P^{(k,0,0,0)}(t,x) &=& k! +\sum_{n>k} R^{(k,0,0,0)}_n(x) \frac{t^{n-k}}{(n-k)!} \\
&=& k! +\sum_{n>k} (kR^{(k,0,0,0)}_{n-1}(x) + (n-k)xR^{(k,0,0,0)}_{n-1}(x)) \frac{t^{n-k}}{(n-k)!}\\
&=& k! + tx \sum_{n>k} R^{(k,0,0,0)}_{n-1}(x) \frac{t^{n-k-1}}{(n-k-1)!}+ 
k \sum_{n > k} R^{(k,0,0,0)}_{n-1}(x) \frac{t^{n-k}}{(n-k)!}.
\end{eqnarray*}
Hence  
\begin{equation}\label{eq:k4}
(1-tx)P^{(k,0,0,0)}(t,x) = k! + 
k \sum_{n > k} R^{(k,0,0,0)}_{n-1}(x) \frac{t^{n-k}}{(n-k)!}.
\end{equation}
Taking the derivative of both sides with respect to $t$, we see 
that
\begin{equation*}
 \frac{\partial}{\partial t}((1-tx)P^{(k,0,0,0)}(t,x))
= k P^{(k,0,0,0)}(t,x) = \frac{k}{1-tx} (1-tx)P^{(k,0,0,0)}(t,x).
\end{equation*}
Thus
\begin{equation}\label{eq:k4-}
\frac{\frac{\partial}{\partial t}((1-tx)P^{(k,0,0,0)}(t,x))}{(1-tx))P^{(k,0,0,0)}(t,x)} =  \frac{k}{1-tx}.
\end{equation}
It then follows that 
\begin{equation*}
\ln((1-tx)P^{(k,0,0,0)}(t,x)) = \frac{-k}{x} \ln(1-tx) + c.
\end{equation*}
Hence, using the fact that $P^{(k,0,0,0)}(0,x) =k!$, we obtain 
that 
\begin{equation}\label{eq:k000}
P^{(k,0,0,0)}(t,x) = k! (1-tx)^{-(\frac{k}{x}+1)}.
\end{equation}

Next, we consider a modification of the notion of a marked mesh pattern.
That is, given a permutation $\sg = \sg_1 \ldots \sg_n \in S_n$, we 
say that $\sg_i$ matches the marked mesh pattern $MMP(=k,0,0,0)$ in 
$\sg$ if and only if, relative to the coordinate system with 
origin at $(i,\sg_i)$, there are exactly $k$ points in $G(\sg)$ 
in quadrant I.  Let $mmp^{(=k,0,0,0)}(\sg)$ denote the number 
of $i$ such that $\sg_i$ matches $MMP(=k,0,0,0)$ and 
$$R_n^{(=k,0,0,0)}(x) = \sum_{\sg \in S_n}x^{mmp^{(=k,0,0,0)}(\sg)}.$$
Then we can use the same reasoning to show 
that 
$$R_n^{(=k,0,0,0)}(x) =n!$$
if $n \leq k$, and, for $n \geq k$, 
\begin{equation}\label{eq:k=2}
R^{(=k,0,0,0)}_{n+1}(x) = (n + x)R^{(=k,0,0,0)}_n(x).
\end{equation}
The difference in this case is that insertion of 1 will 
contribute to $mmp^{(=k,0,0,0)}(\sg^{(i)})$ if and only if 
$i = n+1 -k$. 
Iterating this recursion, we see that 
\begin{equation}\label{eq:k=3}
R^{(=k,0,0,0)}_{k+s}(x) = k! \prod_{i=1}^s(k+i-1+x) 
\end{equation}
for all $s \geq 1$. It is then easy to see using the 
same reasoning as above that if we let 
\begin{equation} \label{P=k0000}
P^{(=k,0,0,0)}(t,x) = \sum_{n\geq k} \frac{t^{n-k}}{(n-k)!} 
\sum_{\sg \in S_n} x^{mmp^{(=k,0,0,0)}(\sg)},
\end{equation}
then 
$$P^{(=k,0,0,0)}(t,x)= k!\frac{1}{(1-t)^{x+k}}.$$

It is also easy to see that a point $\sg_i$  matches 
the pattern $MMP(\emptyset,0,0,0)$ in  $\sg = \sg_1 \ldots \sg_n \in S_n$ 
if and only if $\sg_i$ is 
a right-to-left maximum of $\sg$. It is well-known 
that the number of permutations $\sg \in S_n$ such 
that $\RLmax[\sg] =k$ is equal to the signless Stirling number 
of the first kind $c(n,k)$ which is the number of permutations $\tau \in S_n$ such that 
$\tau$ has $k$ cycles. Thus 
\begin{eqnarray*}
R_n^{(\emptyset,0,0,0)}(x) &=& \sum_{\sg \in S_n} x^{\RLmax[\sg]} =\sum_{k=1}^n c(n,k) x^k \\
&=& x(x+1) \cdots (x+n-1).
\end{eqnarray*}
It is well-known that $\sum_{n \geq 0} \frac{t^n}{n!} \sum_{k=1}^n 
c(n,k)x^k = \left( \frac{1}{1-t}\right)^x$. Thus 
\begin{equation}\label{eq:e000}
R^{(\emptyset,0,0,0)}(t,x) = 1+ \sum_{n \geq 1} R_n^{(\emptyset,0,0,0)}(x)\frac{t^n}{n!} = 
  \left( \frac{1}{1-t}\right)^x.
\end{equation}

\section{$R_n^{(a,b,0,0)}(x)$  and some other cases of two non-zero parameters}

Note that it follows from Lemma \ref{sym} that $R_n^{(a,b,0,0)}(x)$ equals
$$R_n^{(b,0,0,a)}(x) = R_n^{(0,0,a,b)}(x)= 
R_n^{(0,a,b,0)}(x)=R_n^{(b,a,0,0)}(x) =$$
$$R_n^{(0,b,a,0)}(x) = R_n^{(0,0,b,a)}(x)=R_n^{(a,0,0,b)}(x).$$

First, we consider the case where $a,b \geq 1$. 
It is easy to see that $R^{(a,b,0,0)}_n(x) = n!$ if 
$n \leq a+b$.  For $n \geq a+b$, it is easy to see that if $\sg=\sg_1\ldots\sg_n\in S_n$, then
\begin{equation}\label{eq:ab1}
mmp^{(a,b,0,0)}(\sg^{(i)}) = \begin{cases} 
mmp^{(a,b,0,0)}(\sg) & \mbox{if $i \leq b$ or $i > n-a+1$,} \\
1+ mmp^{(a,b,0,0)}(\sg) & \mbox{if $b+1 \leq i \leq  n-a+1$.}
\end{cases}
\end{equation}
That is, clearly the placement of 1 in position $i$ in 
$\sg^{(i)}$ cannot effect the elements 
in quadrant I or II relative to any pair 
$(j,\sg_j+1)$ for 
any $j <i$  or relative to any pair $(j+1,\sg_j+1)$ for any $j \geq i$, 
and the $1$ in position $i$ can contribute  
to $mmp^{(k,0,0,0)}(\sg^{(i)})$ if and only if $b+1 \leq i\leq n-a+1$. 
It then follows that for $n \geq a+b$, 
\begin{equation}\label{eq:ab2}
R^{(a,b,0,0)}_{n+1}(x) = (a+b)R^{(a,b,0,0)}_n(x) + 
(n+1-(a+b))xR^{(a,b,0,0)}_n(x).
\end{equation}
Note this is the same recursion as the recursion 
for $R_{n+1}^{(k,0,0,0)}(x)$ if $a+b =k$. 
Thus it follows that 
\begin{equation}\label{eq:ab3} 
R^{(a,b,0,0)}_{a+b+s}(x) = (a+b)! \prod_{i=1}^s((a+b)+ix) 
\end{equation}
for all $s \geq 1$ and that 
\begin{equation}\label{eq:ab00}
P^{(a,b,0,0)}(t,x) = \sum_{n \geq a+b} \frac{t^{n-a-b}}{(n-a-b)!} 
R_n^{(a,b,0,0)}(x) = 
(a+b)! (1-tx)^{-(\frac{a+b}{x}+1)}
\end{equation}
for all $a,b \geq 1$.

Next we consider the case $R_n^{(1,\emptyset,0,0)}(x)$. 
Clearly, $R_1^{(1,\emptyset,0,0)}(x) =1$.  
For $n \geq 1$ and any $\sg  = \sg_1 \ldots \sg_n \in S_n$, we 
again consider the permutations $\sg^{(i)}$ in $S_{n+1}$. 
For $i=1$, 
$mmp^{(1,\emptyset,0,0)}(\sg^{(1)}) = 1 +
mmp^{(1,\emptyset,0,0)}(\sg)$ and, for $i > 1$,  
$mmp^{(1,\emptyset,0,0)}(\sg^{(1)}) = 
mmp^{(1,\emptyset,0,0)}(\sg)$. It follows that for all $n \geq 1$,
$$R_{n+1}^{(1,\emptyset,0,0)}(x) =(x+n) R_{n}^{(1,\emptyset,0,0)}(x)$$
so that for $n > 1$, 
\begin{equation}\label{eq:1e}
R_n^{(1,\emptyset,0,0)}(x) = (x+1)\cdots (x+n-1).
\end{equation}
It then easily follows that 
\begin{equation}\label{eq3:1e}
R^{(1,\emptyset,0,0)}(t,x) = 1+ \frac{(1-t)^{-x}-1}{x}.
\end{equation}

Next fix $k \geq 2$. Clearly,  
$R_n^{(k,\emptyset,0,0)}(x) =n!$ for $n \leq k$. 
For $n \geq k$ and $\sg \in S_n$, it is easy to see 
that for $i=1$, 
$mmp^{(k,\emptyset,0,0)}(\sg^{(1)}) = 1 +
mmp^{(k,\emptyset,0,0)}(\sg)$ and that for $i > 1$, 
$mmp^{(k,\emptyset,0,0)}(\sg^{(i)}) = 
mmp^{(k,\emptyset,0,0)}(\sg)$. Hence, for all $n \geq k$,
$$R_{n+1}^{(k,\emptyset,0,0)}(x) =(x+n) R_{n}^{(k,\emptyset,0,0)}(x).$$
Thus for $n > k$, 
\begin{equation}\label{eq:1e-}
R_n^{(k,\emptyset,0,0)}(x) = k! (x+k)\cdots (x+n-1).
\end{equation}
It follows that for $k \geq 2$, 
\begin{equation}
P^{(k,\emptyset,0,0)}(t,x) = \sum_{n \geq k} R_n^{(k,\emptyset,0,0)}(x)
\frac{t^{n-k}}{(n-k)!} = k!(1-t)^{-(x+k)} 
\end{equation}
and that 
\begin{multline}\label{kemp00}
R^{(k,\emptyset,0,0)}(t,x) = 1+ \sum_{n \geq 1} R_n^{(k,\emptyset,0,0)}(x)
\frac{t^n}{n!} = \\
\sum_{j=0}^k t^j + \frac{k!}{\prod_{i=1}^{k-1} (x+i)}
\left(R^{(1,\emptyset,0,0)}(t,x) -1-t -\sum_{j=2}^k \frac{t^j}{j!} 
\prod_{i=1}^{j-1} (x+i)\right).
\end{multline}

\section{$q$-analogues to marked mesh patterns considered above}

We let  
\begin{eqnarray*}
\  [n]_{p,q} &=& p^{n-1}+ p^{n-2}q + \cdots + pq^{n-2}+ q^{n-1} = 
\frac{p^n-q^n}{p-q},\\
\ [n]_{p,q}! &=& [1]_q [2]_q \cdots [n]_q, \ \mbox{and} \\
\ \qbin{n}{k}{p,q} &=& \frac{[n]_{p,q} !}{[k]_{p,q}![n-k]_{p,q}!}
\end{eqnarray*}
denote the usual $p,q$-analogues of $n$, $n!$, and $\binom{n}{k}$. 
We shall use the standard conventions that $[0]_{p,q} = 0$ and $[0]_{p,q}! =1$.
Setting $p=1$ in $[n]_{p,q}$, $[n]_{p,q}!$, and 
$\qbin{n}{k}{p,q}$ yields $[n]_q$, $[n]_q!$, and 
$\qbin{n}{k}{q}$, respectively.  
For any permutation $\sg = \sg_1 \ldots \sg_n \in S_n$, we let 
$\inv[\sg]$ equal the number of $1 \leq i < j \leq n$ such that 
$\sg_i > \sg_j$ and $\coinv[\sg]$ equal the number of $1 \leq i < j \leq n$ 
such that  $\sg_i < \sg_j$.

Let 
\begin{equation} \label{qanal1}
R_n^{(a,b,c,d)}(x,q) = \sum_{\sg \in S_n} x^{mmp^{(a,b,c,d)}(\sg)} q^{\coinv}.
\end{equation} 
It turns out that we can easily obtain $q$-analogues 
of the recursions (\ref{eq:k2}) and (\ref{eq:ab2}). 
That is, for any permutation $\sg = \sg_1 \ldots \sg_n \in S_n$, 
it is easy to see that 
\begin{equation}\label{qanal2}
q^{\coinvi} x^{mmp^{(k,0,0,0)}(\sg^{(i)})} = \begin{cases} 
q^{n+1 -i} q^{\coinv} x^{mmp^{(k,0,0,0)}(\sg)} & \mbox{if $i > n-k$}, \\
q^{n+1 -i} q^{\coinv} x^{1+ mmp^{(k,0,0,0)}(\sg)} & \mbox{if $i \leq  n-k$}.
\end{cases}
\end{equation}
It then follows that for $n \geq k$, 
\begin{equation}\label{qanal3}
R^{(k,0,0,0)}_{n+1}(x,q) = [k]_qR^{(k,0,0,0)}_n(x,q) + xq^k[n+1-k]_q
R^{(k,0,0,0)}_n(x,q).
\end{equation}
Thus for $n \leq k$, 
\begin{equation}
R^{(k,0,0,0)}_{n+1}(x,q) =[n]_q!
\end{equation}
and for $s \geq 1$, 
\begin{equation}
R^{(k,0,0,0)}_{k+s}(x,q) =[k]_q! \prod_{i=1}^s ([k]_q +x q^k [i]_q).
\end{equation}

Similarly, we can find a $q$-analogue of the recursion 
(\ref{eq:ab1}).   That is, for any $\sg = \sg_1 \ldots \sg_n \in S_n$, 
it is easy to see that 
\begin{equation}\label{qanal4}
q^{\coinvi} x^{mmp^{(a,b,0,0)}(\sg^{(i)})} = \begin{cases} 
q^{n+1-i} q^{\coinv} x^{mmp^{(a,b,0,0)}(\sg)} & \mbox{if $i \leq b$ or $i > 
n-a+1$}, \\
q^{n+1-i} q^{\coinv} x^{1+ mmp^{(a,b,0,0)}(\sg)} & \mbox{if $b+1 \leq i \leq  n-a+1$}.
\end{cases}
\end{equation}
It then follows that for $n\geq a+b$,
\begin{equation}\label{qanal5}
R^{(a,b,0,0)}_{n+1}(x,q) = ([a]_q+q^{n-b}[b]_q)R^{(k,0,0,0)}_n(x,q) + 
q^a[n+1-(a+b)]_qxR^{(k,0,0,0)}_n(x,q).
\end{equation}
Thus for $n \leq a+b$, 
\begin{equation}
R^{(a,b,0,0)}_{n+1}(x,q) =[n]_q!
\end{equation}
and for $s \geq 1$,
\begin{equation} 
R^{(a,b,0,0)}_{a+b+s}(x,q) =[k]_q!\prod_{i=1}^s ([a]_q +q^{a+i}[b]_q + 
q^ax [i]_q).
\end{equation}
Note that if $a +b =k$ where $a,b \geq 1$  and $n > k$, then  
$R^{(a,b,0,0)}_{n}(x,q) \neq R^{(k,0,0,0)}_{n}(x,q)$.

\section{$R_n^{(k \leq \max,\emptyset,0,0)}(x)$ -- a modification of the notion of a marked mesh pattern}

We say that $\sg_i$ matches the marked  mesh pattern 
$MMP(k\leq \max,\emptyset,0,0)$ in $\sg =\sg_1 \ldots \sg_n \in S_n$  if 
in $G(\sg)$ relative to the coordinate system with origin  
$(i,\sg_i)$, there 
are no points in quadrant II in $G(\sg)$  and 
there are at least $k-1$ points that lie in quadrant I   
to the left of the largest value occurring in quadrant I. Said 
another way, $\sg_i$ matches the pattern 
$MMP(k \leq \max,\emptyset,0,0)$ in $\sg$  if none of $\sg_1, \ldots, \sg_{i-1}$ 
is greater than $\sg_i$ and if $\sg_j$ is the maximum of 
$\sg_{i+1}, \ldots, \sg_n$, then there are at least $k$ points 
among $\sg_{i+1}, \ldots, \sg_j$ which are greater than $\sg_i$.
For $\sg \in S_n$, we let $mmp^{(k\leq \max,\emptyset,0,0)}(\sg)$ 
be the number of $i$ such that $\sg_i$ matches $MMP(k \leq \max,\emptyset,0,0)$ in $\sg$ and we let 
$$R_n^{(k \leq \max,\emptyset,0,0)}(x)= \sum_{\sg \in S_n} x^{mmp^{(k\leq \max,\emptyset,0,0)}(\sg)}.$$

Clearly, $R_n^{(1 \leq \max,\emptyset,0,0)}(x) = R_n^{(1,\emptyset,0,0)}(x)$.
We can compute $R_n^{(k \leq \max,\emptyset,0,0)}(x)$ for 
$k \geq 2$ as follows.
Given a permutation $\sg = \sg_1 \ldots \sg_n \in S_n$, 
let $\sg^{[i]}$ denote the permutation of $S_{n+1}$ that results 
by inserting $n+1$ immediately before $\sg_i$ if $1 \leq i \leq n$ and 
inserting $n+1$ at the end of $\sg$ if $i =n+1$. It is easy to 
see that the set of all $\sg^{[i]}$ for $\sg \in S_n$ contributes 
$\binom{n}{i-1}R_{i-1}^{(k-1,\emptyset,0,0)}(x) (n+1-i)!$ to 
$R_{n+1}^{(k \leq \max,\emptyset,0,0)}(x)$. That is, 
the presence of $n+1$ in the $i$-th position of $\sg^{[i]}$ means 
that none of the elements to the right of $n+1$ in 
$\sg^{[i]}$ can contribute to $mmp^{(k \leq \max,\emptyset,0,0)}(\sg^{[i]})$ 
while an element to the left of $n+1$ in 
$\sg^{[i]}$ will contribute 1 to $mmp^{(k \leq \max,\emptyset,0,0)}(\sg^{[i]})$ if and 
only if it contributes 1 to $mmp^{(k-1,\emptyset,0,0)}(\sg_1 \ldots \sg_{i-1})$. 
It follows that for all $n \geq 0$, 
\begin{equation}\label{eq:1e--}
R_{n+1}^{(k \leq \max,\emptyset,0,0)}(x) = \sum_{i=1}^{n+1} 
(n+1-i)! \binom{n}{i-1}R_{i-1}^{(k-1,\emptyset,0,0)}(x)
\end{equation}
 or,  
equivalently,
\begin{equation}\label{eq2:1e}
\frac{R_{n+1}^{(k \leq \max,\emptyset,0,0)}(x)}{n!}  = \sum_{i=1}^{n+1} 
\frac{R_{i-1}^{(k-1,\emptyset,0,0)}(x)}{(i-1)!}.
\end{equation}
Thus multiplying (\ref{eq2:1e}) by $t^n$ and summing it over all $n \geq 0$, we obtain 
that 
\begin{eqnarray*}
\frac{\partial}{\partial t} R^{(k \leq \max,\emptyset,0,0)}(t,x) &=&
\sum_{n \geq 0} \frac{R_{n+1}^{(k \leq \max,\emptyset,0,0)}(x)t^n}{n!}\\
&=& \sum_{n \geq 0} t^n \sum_{i=1}^{n+1} 
\frac{R_{i-1}^{(k-1,\emptyset,0,0)}(x)}{(i-1)!} \\
&=& \frac{1}{1-t} R^{(k-1,\emptyset,0,0)}(t,x).
\end{eqnarray*}
Hence we obtain the recursion  
\begin{equation}\label{eq2:ke}
R^{(k \leq \max,\emptyset,0,0)}(t,x) = 1+ \int_0^t \frac{1}{1-z}R^{(k-1, \emptyset,0,0)}(z,x)dz.
\end{equation}
Note that one can find an explicit formula for 
$R^{(k-1, \emptyset,0,0)}(z,x)$ by 
using formulas  (\ref{eq3:1e}) and (\ref{kemp00}). 
For example, one can use Mathematica to compute that 
\begin{eqnarray*}
&&R^{(2 \leq \max,\emptyset,0,0)}(t,x) = 1+t+t^2+\frac{1}{6} (5+x) t^3+\\
&&\frac{1}{24} \left(17+6 x+x^2\right) t^4+\frac{1}{120} \left(74+35 x+10 x^2+x^3\right) t^5+\\
&&\frac{1}{720} \left(394+225 x+85 x^2+15 x^3+x^4\right) t^6+\\
&&\frac{1}{5040}\left(2484+1624 x+735 x^2+175 x^3+21 x^4+x^5\right) t^7 +\\
&&\frac{1}{40320}\left(18108+13132 x+6769 x^2+1960 x^3+322 x^4+28 x^5+x^6\right) t^8+\\
&&\frac{1}{362880}\left(149904+118124 x+67284 x^2+22449 x^3+4536 x^4+546 x^5+36 x^6+x^7\right) t^9+ \cdots
\end{eqnarray*}
and 
\begin{eqnarray*}
&&R^{(3 \leq \max,\emptyset,0,0)}(t,x) =
1+t+t^2+t^3+\frac{1}{24} 2(11+x) t^4+\frac{1}{120} 2\left(50+9 x+x^2\right) t^5+\\
&&\frac{1}{720} 2\left(274+71x+14 x^2+x^3\right) t^6+\frac{1}{5040} 2(1764+580x+ 155 x^2+20 x^3+x^4) t^7+\\
&&\frac{1}{40320} 2(13068+5104 x+1665x^2+295 x^3+27x^4+x^5) t^8+\\
&&\frac{1}{362880} 2(100584+48860 x+18424 x^2+4025x^3+511 x^4+35x^5+x^6) t^9+ \cdots.
\end{eqnarray*}
and 
\begin{eqnarray*}
&&R^{(4 \leq \max,\emptyset,0,0)}(t,x) =
1+t+t^2+t^3+t^4+\frac{1}{120} 6\left(19+x\right) t^5+\frac{1}{720} 6\left(107+12x+x^2\right) t^6+\\
&&\frac{1}{5040} 6(702+119x+18x^2+x^3) t^7+\frac{1}{40320} 6(5274+1175x+245x^2+25x^3+x^4) t^8+\\
&&\frac{1}{362880} 6(44712+12154x+3135x^2+445x^3+33x^4+x^5) t^9+ \cdots.
\end{eqnarray*}

There are several of the coefficients of 
$R^{(k \leq \max,\emptyset,0,0)}_n(x)$ that we can explain. 
First, we claim that  
\begin{equation}
R^{(k \leq \max,\emptyset,0,0)}_n(x)|_{x^{n-k}} = (k-1)! \ \mbox{for} \ n \geq k+1.
\end{equation}
That is, it is easy to see that for $n \geq k+1$, the only permutations 
$\sg = \sg_1 \ldots \sg_n \in S_n$ such that $mmp^{(k \leq \max,\emptyset,0,0)}(\sg) =n-k$ 
 must have 
$\sg_i =i$ for $i =1, \ldots, n-k$, $\sg_n =n$ and 
$\sg_{n-k+1} \ldots \sg_{n-1}$ be some permutation 
of $(n-k+1), (n-k+2), \ldots, (n-1)$. 

Next, we claim that  
\begin{equation}
R^{(k \leq \max,\emptyset,0,0)}_n(x)|_{x^{n-k-1}} = (k-1)! \left( \binom{n}{2} -\binom{k-1}{2}\right)\ \mbox{for} \ n \geq k+2.
\end{equation}
First observe that any permutation $\sg = \sg_1 \ldots \sg_n \in S_n$ 
such that $n \in \{\sg_1, \ldots, \sg_{n-2}\}$ cannot 
have $mmp^{(k \leq \max,\emptyset,0,0)}(\sg) =n-k-1$. 
 Now assume that $n \geq k+2$ and   
$\sg = \sg_1 \ldots \sg_n \in S_n$ is such that 
$mmp^{(k \leq \max,\emptyset,0,0)}(\sg) =n-k-1$. 
Thus it must be that case that $n = \sg_{n-1}$ or $n = \sg_n$.
If $\sg_{n-1} =n$, then it must be the case that 
$mmp^{(k \leq \max,\emptyset,0,0)}(\red[\sg_1 \ldots \sg_{n-1}]) = n-k-1$. 
In this case, we have $n-1$ choices for $\sg_n$ and 
$(k-1)!$ choices for $\red[\sg_1 \ldots \sg_{n-1}]$ by our 
argument above. Thus there are $(n-1)((k-1)!)$ such elements. 
Next assume that $n = \sg_n$.  Then in this case, the only 
way to produce a  
$\sg = \sg_1 \ldots \sg_n \in S_n$ such that 
$mmp^{(k \leq \max,\emptyset,0,0)}(\sg) =n-k-1$ is  by  
starting with a permutation 
$$\tau = 1~2~ \ldots (n-k) \tau_{n-k+1} \ldots \tau_{n-1}~n$$
 where 
$ \tau_{n-k+1} \ldots \tau_{n-1}$ is some permutation 
of $(n-k+1), (n-k+2), \ldots, (n-1)$ so that 
$mmp^{(k \leq \max,\emptyset,0,0)}(\tau) =n-k$ and then 
take some $i \in \{1, \ldots, n-k\}$ in $\tau$ and move 
$i$ immediately after one of 
$\tau_{i+1}, \ldots, \tau_{n-1}$. This will ensure that 
$i$ does not match the pattern $MMP(k\leq \max,\emptyset,0,0)$ in 
the resulting permutation because 
$i$ will have a larger element to its left. 
Each such 
$\tau$ gives rise of $(n-2) +(n-3) + \cdots + (k-1) = \binom{n-1}{2} - \binom{k-1}{2}$ such  permutations.  
Thus the number of $\sg \in S_n$ with $\sg_n =n$ and 
$mmp^{(k \leq \max,\emptyset,0,0)}(\sg) =n-k-1$ is 
$(k-1)!\left(\binom{n-1}{2} - \binom{k-1}{2}\right)$. 
It thus follows that the number of $\sg \in S_n$ with  
$mmp^{(k \leq \max,\emptyset,0,0)}(\sg) =n-k-1$ is 
$(k-1)!\left(\binom{n}{2} - \binom{k-1}{2}\right)$.

Next we observe that the sequence 
$(R_n^{(2 \leq \max,\emptyset,0,0)}(0))_{n \geq 1}$ is 
$$1,2,5,17,74,394,2484,18108,149904, \ldots.$$
This is sequence A000774 from the OEIS where the $n$-th 
term in the sequence is $n!(1+\sum_{i=1}^n \frac{1}{i})$.  
The sequence $(R_n^{(2 \leq \max,\emptyset,0,0)}(x)|_x)_{n \geq 3}$ is 
$$1,6,35,225,1624,13,132,118124 \ldots$$
which is sequence A000399 in the OEIS.  The $n$-term of this 
sequence is $c(n,3)$ which is 
the number of permutations $\sg \in S_n$ which has 
three cycles. We can give a direct proof of this result. That is, 
suppose that $n \geq 3$ and $\sg$ is a permutation of $S_n$ 
with 3 cycles $C_1,C_2,C_3$ where have arranged 
the cycles  so that 
the largest element in each cycle is on the left and 
we order the cycles by increasing largest elements.
Then we let $\overline{\sg}$ be the result of erasing the parentheses 
and commas in $C_1 C_2 C_3$.  
For example, if $\sg= (5,3,4,1)~(8,2,6)~(9,7)$, then 
$\overline{\sg} = 5~3~4~1~8~2~6~9~7$. 
It is easy to see that under this map only the first element 
of $\overline{\sg}$ matches the pattern 
$MMP(2 \leq \max,\emptyset,0,0)$ in $\overline{\sg}$.  

We claim that every $\tau 
\in S_n$ such that $mmp^{(2 \leq \max,\emptyset,0,0)}(\tau) = 1$ is 
equal to $\overline{\sg}$ for some $\sg \in S_n$ that has 
two three cycles. That is, suppose that  $mmp^{(2 \leq \max,\emptyset,0,0)}(\tau) = 1$ where $\tau = \tau_1 \ldots \tau_n \in S_n$  and 
$\tau_i$ matches the pattern $MMP(2 \leq \max,\emptyset,0,0)$ in $\tau$. 
First we claim that $i=1$. It cannot be that 
any of $\tau_1, \ldots, \tau_{i-1}$ are greater than $\tau_i$ since 
otherwise 
$\tau_i$ would  not match the pattern $MMP(2 \leq \max,\emptyset,0,0)$ in 
$\tau$. 
But if $i \neq 1$, then $\tau_1 < \tau_i$ and $\tau_1$ would 
match the pattern $MMP(2 \leq \max,\emptyset,0,0)$ in $\tau$ which would 
mean that  $mmp^{(2 \leq \max,\emptyset,0,0)}(\tau) \geq 2$.
Next suppose that $\tau_k =n$ and $\tau_j$ is the maximum 
element of $\{\tau_2, \ldots, \tau_{k-1}\}$. It must be 
that case that $\tau_1 < \tau_j$ since otherwise $\tau_1$ would not match the pattern $MMP(2 \leq \max,\emptyset,0,0)$ in $\tau$. We claim that 
$\tau_2, \ldots, \tau_{j-1}$ must all be less than $\tau_j$. That is, if 
$\tau_r$ is the maximum of $\{\tau_2, \ldots, \tau_{j-1}\}$ and 
$\tau_r > \tau_1$, then we would have $\tau_r < \tau_j < \tau_k =n$ 
so that $\tau_r$ would match the pattern $MMP(2 \leq \max,\emptyset,0,0)$ 
in $\tau$ which would 
imply that  $mmp^{(2 \leq \max,\emptyset,0,0)}(\tau) \geq 2$. It 
then follows that 
if $\sg = (\tau_1, \ldots, \tau_{j-1})~(\tau_j, \ldots, \tau_{j-1})~(\tau_k, \ldots, \tau_n)$, then $\overline{\sg} = \tau$.  Thus we have proved 
that for 
$n \geq 3$, $R_n^{(2 \leq \max,\emptyset,0,0)}(x)|_x =c(n,3)$.

The sequence $(\frac{1}{2}R_n^{(3 \leq \max,\emptyset,0,0)}(0))_{n \geq 2}$ is 
$$1,3,11,50,274,1764,13068,109584,1026576,  \ldots$$ which 
is sequence A000254 in the OEIS whose $n$-th term is $c(n+1,2)$ which 
is the number of permutations of $S_{n+1}$ with 2 cycles. 
We shall give a direct proof  
of the fact that $R_n^{(3 \leq \max,\emptyset,0,0)}(0) = 2c(n,2)$ for $n \geq 2$. 
That is, suppose that $\sg = \sg_1 \ldots \sg_n$ is a permutation 
in $S_n$ such that 
$mmp^{(3 \leq \max,\emptyset,0,0)}(\sg) =0$.  Let 
$a(\sg)$ be the $i$ such that $\sg_i=n$.  If $a(\sg) \neq 1$, then 
let $b(\sg) =j$ where $\sg_j = \max(\{\sg_1, \ldots, \sg_{a(\sg)-1}\})$. 
If $b(\sg) \neq 1$, then let $c(\sg) =i$ where $\sg_i = 
\max(\{\sg_1, \ldots, \sg_{b(\sg)-1}\})$. 
First we claim that if $c(\sg)$ is defined, then $c(\sg)=1$.  That is, 
if $c(\sg) \neq 1$, then $\sg_1 < \sg_{c(\sg)} < \sg_{b(\sg)} < 
\sg_{c(\sg)} =n$ so that $\sg_1$ would match the pattern $MMP(3 \leq \max,\emptyset,0,0)$.  We also claim that if $c(\sg)$ is defined, then 
$\sg_1$ must be the second smallest element among 
$\{\sg_1, \ldots, \sg_{a(\sg)-1}\}$ otherwise 
$\sg_1$ would match the pattern $MMP(3 \leq \max,\emptyset,0,0)$ in $\sg$.  
Thus we have three possible cases for a $\sg \in S_n$ such that 
$mmp^{(3 \leq \max,\emptyset,0,0)}(\sg) =0$.\\
\ \\
{\bf Case 1.} $a(\sg) = 1$.\\
There are clearly $(n-1)!$ such permutations as $\sg_2, \ldots, \sg_n$ 
can be any arrangement of $1, \ldots, n-1$. \\
\ \\
{\bf Case 2.} $b(\sg) =1$.\\
\ \\
{\bf Case 3.} $c(\sg) =1$ and $\sg_1$ is the second 
smallest element in $\{\sg_1, \ldots, \sg_{a(\sg) -1}\}$.

Next we define a map $\theta$ which takes the permutations 
in Cases 2 and 3 into the set of permutations $\tau$ of $S_n$ which 
have 2 cycles $C_1C_2$ where we have arranged the cycles so 
that the largest element in each cycle is on the left and 
we have $\max(C_1) < \max(C_2) =n$. 
In Case 2, we let 
$$
\theta(\sg) = (\sg_{b(\sg)}, \ldots, \sg_{a(\sg) -1})~ (\sg_{a(\sg)}, \ldots , \sg_n).$$
In Case 3, we let 
$$
\theta(\sg) = 
(\sg_{b(\sg)}, \ldots, \sg_{a(\sg) -1}, \sg_{c(\sg)}, \ldots, 
\sg_{b(\sg)-1})~ (\sg_{a(\sg)}, \ldots , \sg_n).$$

Now suppose that $\tau$ is a permutation with 
two cycles $C_1C_2$, $C_1 = (\alpha_1, \ldots, \alpha_k)$ and $ C_2=(\beta_1, \ldots, \beta_{n-k})$,  where  $|C_1| = k > 1$. 
Again we assume 
that we have  arranged the cycles so that the largest element 
in the cycle is on the left and  $\max(C_1) < \max(C_2) =n$. 
Then  there exists a $\sg$ in 
Case 2 such that $\theta(\sg) = \tau$ and a $\gamma$ in Case 3 such 
that $\theta(\gamma) = \tau$. That is, if $\alpha_j$ is the 
second smallest element in $\{\alpha_1, \ldots, \alpha_k\}$, 
then 
$$\sg = \alpha_1, \ldots, \alpha_k \beta_1 \ldots \beta_{n-k}$$ is 
a permutation in Case 2 such that $\theta(\sg)=\tau$ and 
$$\gamma = \alpha_j \ldots \alpha_k \alpha_1, \ldots, \alpha_{j-1} \beta_1 \ldots \beta_{n-k}$$
is a permutation in Case 3 such that $\theta(\gamma) = \tau$. 
The only permutations with 2 cycles $C_1C_2$ that we have not accounted 
for are the permutations $\tau$ where $|C_1| =1$. But if $|C_1| =1$, 
then $\sg = \alpha_1\beta_1 \ldots \beta_{n-1}$ is a permutation 
in Case 2 such that $\theta(\sg) =\tau$. It follows 
that $\theta$ shows that the number of permutations 
in Cases 2 and 3 is equal to $2c(n,2)$ minus the number of 
$\tau$ such that $|C_1| =1$.  However, it is easy to see 
that the number of $\tau$ with two cycles $C_1C_2$ such that $|C_1| =1$ is 
$(n-1)!$ which is equal to the number of permutations in Case 1. 
Thus we have shown that $R_n^{(3 \leq \max,\emptyset,0,0)}(0) = 2c(n,2)$.

The sequence 
$(\frac{1}{6}R_n^{(\leq 4 \max,\emptyset,0,0)}(0))_{n \geq 3}$ which is 
$$1,4,19,107,702,5274,44712,422,568 \ldots $$
does not appear in the OEIS. However, the sequence $(\frac{1}{6}
R_n^{(4 \leq \max,\emptyset,0,0)}(x)|_x)_{n \geq 5}$ which is 
$$1,12,119,1175,12154,133938, \ldots $$
seems to be A001712 in the OEIS whose $n$-th term is $\sum_{k=0}^n 
(-1)^{n+k} \binom{k+2}{2}3^k s(n+2,k+2)$ where 
$s(n,k)$ is the Stirling number of the first kind. 

\begin{problem} Can we prove this formula (directly)? \end{problem}

\section{$R_n^{(a,0,b,0)}(x)=R_n^{(b,0,a,0)}(x)=R_n^{(0,a,0,b)}(x)=R_n^{(0,b,0,a)}(x)$}

The equalities in the section title are true by Lemma \ref{sym}. We will consider $R_n^{(a,0,b,0)}(x)$.

In this case, we do not know how to compute the generating 
function $R^{(1,0,1,0)}(t,x)$. 
However, we can develop a recursion 
to compute $R_n^{(1,0,1,0)}(x)$.  That is, 
for any permutation $\sg = \sg_1 \ldots \sg_n \in S_n$, 
let 
\begin{eqnarray*}
A_i(\sg) &=& \chi(\sg_i \ \mbox{matches the pattern} \ MMP(1,0,1,0)\ \mbox{in} \ \sg) \ 
\mbox{and} \\
B_i(\sg) &=& \chi(\sg_i \ \mbox{matches the pattern} \ 
MMP(\emptyset,0,1,0) \ \mbox{in} \ \sg)
\end{eqnarray*}
where for any statement $A$, $\chi(A) =1$ if $A$ is true and 
$\chi(A) =0$ if $A$ is false. 
Then let 
\begin{equation}
F^{(1,0,1,0)}_n(x_2, \ldots, x_n;y_2, \ldots, y_n) = 
\sum_{\sg \in S_n} \prod_{i=2}^n x_i^{A_i(\sg)} y_i^{B_i(\sg)}.
\end{equation}
Note that $(1,\sg_1)$ never contributes to $mmp^{(1,0,1,0)}(\sg)$ or 
$mmp^{(\emptyset,0,1,0)}(\sg)$ so that \\
$F^{(1,0,1,0)}_n(x_2, \ldots, x_n;y_2, \ldots, y_n)$ records all 
the information about the contributions of \\
$\sg_2, \ldots, \sg_n$ to $mmp^{(1,0,1,0)}(\sg)$ or 
$mmp^{(\emptyset,0,1,0)}(\sg)$ as $\sg$ ranges over $S_n$. 

Recall that given a permutation $\sg = \sg_1 \ldots \sg_n \in S_n$, 
$\sg^{[i]}$ denotes the permutation of $S_{n+1}$ that results 
by inserting $n+1$ immediately before $\sg_i$ if $1 \leq i \leq n$ and 
inserting $n+1$ at the end of $\sg$ if $i =n+1$.
First consider all the permutations  $\sg^{[1]}$ as $\sg$ ranges over 
$S_n$. Clearly, since each of these permutations start with $n+1$, 
the first element of $\sg^{[1]}$ which is $n+1$ 
 does not contribute to either $mmp^{(1,0,1,0)}(\sg^{[1]})$ or 
$mmp^{(\emptyset,0,1,0)}(\sg^{[1]})$. Moreover $n+1$ does not 
effect whether any other $\sg_i$ contributes to 
$mmp^{(1,0,1,0)}(\sg^{[1]})$ or 
$mmp^{(\emptyset,0,1,0)}(\sg^{[1]})$ but it does shift the 
corresponding indices by 1. Thus the second element 
of $\sg^{[1]}$ which is $\sg_1$ does not contribute to either 
$mmp^{(1,0,1,0)}(\sg^{[1]})$ or 
$mmp^{(\emptyset,0,1,0)}(\sg^{[1]})$. However, for the $(i+1)$-st element 
$\sg^{[1]}$ which is the old $\sg_i$, we have 
\begin{enumerate}
\item $A_{i+1}(\sg^{[1]}) =1$ if and only if $A_{i}(\sg) =1$ and 
\item $B_{i+1}(\sg^{[1]}) =1$ if and only if $B_{i}(\sg) =1$.
\end{enumerate}
Thus it follows that the contribution 
of the permutations of the form $\sg^{[1]}$ to \\ 
$F^{(1,0,1,0)}_{n+1}(x_2,\ldots,x_{n+1};y_2,\ldots, y_{n+1})$ is just 
$F^{(1,0,1,0)}_n(x_3, \ldots, x_{n+1};y_3, \ldots, y_{n+1})$. 
For $i \geq 2$, again 
consider all the permutations  $\sg^{[i]}$ as $\sg$ ranges over 
$S_n$. For $j < i$, if $A_j(\sg) =1$, then 
$A_j(\sg^{[i]})  =1$ and if $B_j(\sg) =1$, then 
$A_j(\sg^{[i]})  =1$.  For $j \geq i$, if $A_j(\sg) =1$, then 
$A_{j+1}(\sg^{[i]})  =1$  and if $B_j(\sg) =1$, then 
$B_{j+1}(\sg^{[i]})  =1$. Moreover the fact that 
$n+1$ is in position $i$ in $\sg^{[i]}$ means that 
$B_i(\sg^{[i]})  =1$. Since at most one 
of $A_i(\sg)$ and $B_i(\sg)$ can equal 1, it follows that 
the contribution 
of the permutations of the form $\sg^{[i]}$ to 
$F_{n+1}(x_2,\ldots,x_{n+1};y_2,\ldots, y_{n+1})$ is just 
$$y_iF^{(1,0,1,0)}_n(x_2, \ldots,x_{i-1},x_{i+1},\ldots,  x_{n+1};x_2,\ldots,x_{i-1},y_{i+1}, \ldots, y_{n+1}).$$ Thus 
\begin{eqnarray}\label{1010rec}
&&F^{(1,0,1,0)}_{n+1}(x_2,\ldots, x_{n+1};y_2, \ldots, y_{n+1}) =  F^{(1,0,1,0)}_n(x_3, \ldots, x_{n+1};y_3, \ldots, y_{n+1}) + \nonumber \\
&&  \sum_{i=2}^{n+1} y_iF^{(1,0,1,0)}_n(x_2, \ldots,x_{i-1},x_{i+1},\ldots,  x_{n+1};x_2,\ldots,x_{i-1},y_{i+1}, \ldots, y_{n+1}).
\end{eqnarray}
It is then easy to see that 
\begin{equation}\label{1010rec2}
R_n^{(1,0,1,0)}(x) = F^{(1,0,1,0)}_n(x, \ldots, x;1, \ldots,1).
\end{equation}
Using (\ref{1010rec}) and (\ref{1010rec2}), one can compute 
that \\
$R_1^{(1,0,1,0)}(x) =1$,\\
$R_2^{(1,0,1,0)}(x) =2$,\\
$R_3^{(1,0,1,0)}(x) =5+x$,\\
$R_4^{(1,0,1,0)}(x) =14+8x+2x^2$,\\
$R_5^{(1,0,1,0)}(x) =42+46x+26x^2+6x^3$,\\
$R_6^{(1,0,1,0)}(x) =132+232x+220x^2+112x^3+24x^5$,\\
$R_7^{(1,0,1,0)}(x) =429+1093x+1527x^2+1275x^3+596x^4+120x^5$, and \\
$R_8^{(1,0,1,0)}(x) =1430+4944x+9436x^2+11384x^3+8638x^4+3768x^5+720x^6$.\\

It is easy to see that a permutation $\sg$ avoids 
$MMP(1,0,1,0)$ if and only if it avoids the pattern $123$, that is, if $\sg\in S_n(123)$. 
From a well-known fact (see, e.g. \cite{kit}), it follows that 
the sequence $(R_n^{(1,0,1,0)}(0))_{n \geq 1}$ 
is the Catalan numbers. Thus it is not surprising that it 
is difficult to find a simple expression for $R^{(1,0,1,0)}(t,x)$ since 
then  $R^{(1,0,1,0)}(t,0)$ would give us an exponential generating 
function for the Catalan numbers which is not known.  
It is also easy to see that for any $n > 2$, the most occurrences 
of the pattern $MMP(1,0,1,0)$ occurs when $\sg_1 =1$ and $\sg_n=n$. 
Clearly, there are $(n-2)!$ such permutations which explains 
the highest coefficients in the $(R_n^{(1,0,1,0)}(x))_{n \geq 1}$. 
We note that neither the sequence $(R_n^{(1,0,1,0)}(x)|_x)_{n \geq 3}$ 
nor the sequence $(R_n^{(1,0,1,0)}(x)|_{x^{n-3}})_{n \geq 3}$ appear 
in the OEIS.

We can use similar reasoning  to develop a recursion 
to compute $R_n^{(1,0,a,0)}(x)$.  That is, 
for any permutation $\sg = \sg_1 \ldots \sg_n \in S_n$, 
let 
\begin{eqnarray*}
A_i^a(\sg) &=& \chi(\sg_i \ \mbox{matches the pattern} \ MMP(1,0,a,0) \ \mbox{in} \ \sg) \ 
\mbox{and} \\
B_i^a(\sg) &=& \chi(\sg_i \ \mbox{matches the pattern} \ 
MMP(\emptyset,0,a,0) \ \mbox{in} \ \sg).
\end{eqnarray*}
Then let 
\begin{equation}
F_n^{(1,0,a,0)}(x_2, \ldots, x_n;y_2, \ldots, y_n) = 
\sum_{\sg \in S_n} \prod_{i=2}^n x_i^{A_i^a(\sg)} y_i^{B_i^a(\sg)}.
\end{equation}
We can now use the same analysis that we did to prove 
(\ref{1010rec}) to prove that 
\begin{eqnarray}\label{10a0rec}
&&F_{n+1}^{(1,0,a,0)}(x_2,\ldots, x_{n+1};y_2, \ldots, y_{n+1}) 
= F_n^{(1,0,a,0)}(x_3, \ldots, x_{n+1};y_3, \ldots, y_{n+1}) + \nonumber \\
&& \sum_{i=2}^{a} F_n^{(1,0,a,0)}(x_2, \ldots,x_{i-1},x_{i+1},\ldots,  x_{n+1};x_2,\ldots,x_{i-1},y_{i+1}, \ldots, y_{n+1}) + \nonumber \\
&&\sum_{i=a+1}^n  y_i F_n^{(1,0,a,0)}(x_2, \ldots,x_{i-1},x_{i+1},\ldots,  x_{n+1};x_2,\ldots,x_{i-1},y_{i+1}, \ldots, y_{n+1}).
\end{eqnarray}
The difference in this case is that if $n \geq a$, then 
inserting $n+1$ immediately 
before $\sg_1, \ldots, \sg_a$ does not give rise to an extra 
factor of $y_i$, while inserting $n+1$ immediately before 
$\sg_{a+1}, \ldots ,\sg_n$ or at the end does give rise to a factor 
of $y_i$ since then $n+1$ matches the pattern $MMP(\emptyset,0,a,0)$. 
It is then easy to see that 
\begin{equation}\label{10a0rec2}
R_n^{(1,0,a,0)}(x) = F_n^{(1,0,a,0)}(x, \ldots, x;1, \ldots,1).
\end{equation}
Using (\ref{10a0rec}) and (\ref{10a0rec2}) one can compute 
that for $a=2$,\\
$R_1^{(1,0,2,0)}(x) =1$,\\
$R_2^{(1,0,2,0)}(x) =2$,\\
$R_3^{(1,0,2,0)}(x) =6$,\\
$R_4^{(1,0,2,0)}(x) =22+2x$,\\
$R_5^{(1,0,2,0)}(x) =90+26x+4x^2$,\\
$R_6^{(1,0,2,0)}(x) =394+232x+82x^2 +12x^3$,\\
$R_7^{(1,0,2,0)}(x) =1806+1776x+1062x^2+348x^3+48x^4$, and \\
$R_8^{(1,0,2,0)}(x) =8558+12546x+11118x^2+6022x^3+1836x^4+240x^5$.\\

Again, it is easy to see that one obtains the maximum 
number of $MMP(1,0,2,0)$-matches in a permutation 
$\sg = \sg_1 \ldots \sg_n \in S_n$ if 
either $\sg_1 =1,\sg_2 =2$, and $\sg_n =n$ or  $\sg_1 =2,\sg_2 =1$, and $\sg_n =n$. 
It follows that for $n \geq 4$, the coefficient of 
the highest power of $x$ in $R_n^{(1,0,2,0)}(x)$ is 
$2(n-3)!$.  The sequence $(R_n^{(1,0,2,0)}(0))_{n \geq 1}$ whose 
initial terms are 
$$1, 2, 6, 22, 90, 394, 1806, 8558,\ldots$$
is the so-called {\em large Schr\"oder numbers}, A006318 in the OEIS. This follows from the fact that avoiding $MMP(1,0,2,0)$ is the same as avoidance simultaneously of the patterns 1234 and 2134, a known case (see \cite[Table 2.2]{kit}). The same numbers count so-called {\em separable permutations} (those avoiding simultaneously the patterns 2413 and 3142) and other eight non-equivalent modulo trivial bijections classes of avoidance of two (classical) patterns of length 4 (see \cite[Table 2.2]{kit}). We note that neither the sequence $(R_n^{(1,0,2,0)}(x)|_x)_{n \geq 4}$ 
nor the sequence $(R_n^{(1,0,2,0)}(x)|_{x^{n-4}})_{n \geq 4}$ appear 
in the OEIS.

One can also develop a recursion  
to compute $R_n^{(b,0,a,0)}(x)$ when $b \geq 2$ but 
it is more complicated. That is, in  such 
a case, we have to keep track of 
the patterns $MMP(\emptyset,0,a,0)$, $MMP(1,0,a,0), \ldots, 
MMP(b,0,a,0)$.  We will show how this works 
in the case when $a=b =2$. 
That is, 
for any permutation $\sg = \sg_1 \ldots \sg_n \in S_n$, 
let 
\begin{eqnarray*}
A_i^{2,2}(\sg) &=& \chi(\sg_i \ \mbox{matches the pattern} \ MMP(2,0,2,0) \ \mbox{in} \ \sg), \\
B_i^{2,2}(\sg) &=& \chi(\sg_i \ \mbox{matches the pattern} \ 
MMP(1,0,2,0) \ \mbox{in} \ \sg, \ \mbox{but does not match the pattern}  \\
&&MMP(2,0,2,0) \ \mbox{in} \ \sg)\ \ \mbox{and}\\
C_i^{2,2}(\sg) &=& \chi(\sg_i \ \mbox{matches the pattern} \ 
MMP(\emptyset,0,2,0) \ \mbox{in} \ \sg).
\end{eqnarray*}
Then let 
\begin{equation}
G_n^{(2,0,2,0)}(x_3, \ldots, x_n;y_3, \ldots, y_n;z_3, \ldots, z_n) = 
\sum_{\sg \in S_n} \prod_{i=3}^n x_i^{A_i^{2,2}(\sg)} y_i^{B_i^{2,2}(\sg)}
z_i^{C_i^{2,2}(\sg)}.\end{equation}

Note that for any $\sg = \sg_1 \ldots \sg_n \in S_n$, we cannot have either $\sg_1$ or $\sg_2$ match any of the patterns 
$MMP(\emptyset,0,2,0)$, $MMP(1,0,2,0)$, or $MMP(2,0,2,0)$ in $\sg$ so that  
is why we do not need variables with subindices 1 or 2.  
Clearly, for $n < 3$,  
$$G_n^{(2,0,2,0)}(x_3, \ldots, x_n;y_3, \ldots, y_n;z_3, \ldots, z_n) =n!.$$   
For $n=3$, we can never have a pattern match 
for $MMP(2,0,2,0)$ or $MMP(1,0,2,0)$. We can have a pattern match 
for $MMP(\emptyset,0,2,0)$ only if $\sg = 123$ or 
$\sg =213$. Thus $G_3^{2,2}(x_3,y_3,z_3) = 4+2z_3$.

Consider all the permutations  $\sg^{[i]}$ where $i \in \{1,2\}$
as $\sg$ ranges over 
$S_n$. Clearly, since each of these permutations start with $n+1$ or 
has its second element equal to $n+1$, 
$n+1$ does not contribute to either $mmp^{(2,0,2,0)}(\sg^{[i]})$,  
$mmp^{(1,0,2,0)}(\sg^{[i]})$, or 
$mmp^{(\emptyset,0,2,0)}(\sg^{[i]})$. Moreover $n+1$ does not 
effect whether any other $\sg_i$ contributes to 
 $mmp^{(2,0,2,0)}(\sg^{[i]})$,  $mmp^{(1,0,2,0)}(\sg^{[i]})$,  or 
$mmp^{(\emptyset,0,2,0)}(\sg^{[i]})$, but it does shift the index 
over by 1. Thus the contribution 
of the permutations of the form $\sg^{[1]}$ and $\sg^{[2]}$ to  
$G_{n+1}^{(2,0,2,0)}(x_3,\ldots,x_{n+1};y_3,\ldots, y_{n+1};z_3,\ldots, z_{n+1})$ is just 
$2G_{n}^{(2,0,2,0)}(x_4, \ldots, x_{n+1};y_4, \ldots, y_{n+1},z_4, \ldots, z_{n+1})$. 
Fix $i \geq  3$ and  
consider all the permutations  $\sg^{[i]}$ as $\sg$ ranges over 
$S_n$. For $j < i$, if $A_j(\sg) =1$, then 
$A_j(\sg^{[i]})  =1$, if $B_j(\sg) =1$, then 
$A_j(\sg^{[i]})  =1$ and if $C_j(\sg) =1$, then 
$B_j(\sg^{[i]})  =1$.  For $j \geq i$, if $A_j(\sg) =1$, then 
$A_{j+1}(\sg^{[i]})  =1$, if $B_j(\sg) =1$, then 
$B_{j+1}(\sg^{[i]})  =1$ and if $C_j(\sg) =1$, then 
$C_{j+1}(\sg^{[i]})  =1$. Moreover, the fact that 
$n+1$ is in position $i$ in $\sg^{[i]}$ means that 
$C_i(\sg^{[i]})  =1$. Since at most one 
of $A_i(\sg)=1$, $B_i(\sg) =1$, and $C_i(\sg) =1$ can hold for any $i$, 
it follows that 
the contribution 
of the permutations of the form $\sg^{[i]}$ to 
$G_{n+1}^{(2,0,2,0)}(x_3,\ldots,x_{n+1};y_3,\ldots, y_{n+1};z_3,\ldots, z_{n+1})$ is just 
$z_iF_n(x_3, \ldots,x_{i-1},x_{i+1},\ldots,  x_{n+1};x_3,\ldots,x_{i-1},y_{i+1}, \ldots, y_{n+1};y_3,\ldots,y_{i-1},z_{i+1}, \ldots, z_{n+1})$. Thus 
\begin{multline}\label{2020rec}
G_{n+1}^{(2,0,2,0)}(x_3,\ldots, x_{n+1};y_3, \ldots, y_{n+1};z_3, \ldots, z_{n+1}) 
=  \\
2 G_{n+1}^{(2,0,2,0)}(x_4,\ldots, x_{n+1};y_4, \ldots, y_{n+1};z_4, \ldots, z_{n+1}) \ +  \\
\sum_{i=3}^{n+1} z_iG_n^{(2,2)}(x_3, \ldots,x_{i-1},x_{i+1},\ldots,  x_{n+1};x_2,\ldots,x_{i-1},y_{i+1}, \ldots, y_{n+1};y_3,\ldots,y_{i-1},z_{i+1}, \ldots, z_{n+1}). 
\end{multline}
It is then easy to see that 
\begin{equation}\label{2020rec2}
R_n^{(2,0,2,0)}(x) = G_n^{(2,0,2,0)}(x, \ldots, x;1, \ldots,1;1, \ldots,1).
\end{equation}
Using (\ref{2020rec}) and (\ref{2020rec2}) one can compute 
that \\
$R_1^{(2,0,2,0)}(x) =1$,\\
$R_2^{(2,0,2,0)}(x) =2$,\\
$R_3^{(2,0,2,0)}(x) =6$,\\
$R_4^{(2,0,2,0)}(x) =24$,\\
$R_5^{(2,0,2,0)}(x) =116+4x$,\\
$R_6^{(2,0,2,0)}(x) =632+80x+8x^2$,\\
$R_7^{(2,0,2,0)}(x) =3720+1056x+240x^2+24x^3$, and \\
$R_8^{(2,0,2,0)}(x) =23072+11680x+4480x^2+992x^3+96x^4$.\\

Again, it is easy to see that one obtains the maximum 
number of $MMP(2,0,2,0)$-matches in a permutation 
$\sg = \sg_1 \ldots \sg_n \in S_n$ if 
$\{\sg_1,\sg_2\} =\{1,2\}$ and $\{\sg_{n-1},\sg_n\}= \{n-1,n\}$.
It follows that for $n \geq 5$, the coefficient of 
the highest power of $x$ in $R_n^{(1,0,2,0)}(x)$ is 
$4(n-4)!$. In this case, the sequence 
$(R^{(2,0,2,0)}(0))_{n \geq 1}$ does not appear in the OEIS, but 
it clearly counts the number of permutations which avoid 
the permutations $12345$, $21345$, $12354$, and $21354$. 

While we cannot find a generating function for 
$R^{(1,0,1,0)}(t,x)$, we can find a restricted version of 
this generating function. That is, let 
$S_n(1 \rightarrow n)$ denote the set of permutations 
$\sg = \sg_1 \ldots \sg_n$ such that 1 appears to the 
left of  $n$  in $\sg$.  Then let 
\begin{equation}\label{B1010}
B^{(1,0,1,0)}(t,x) = \sum_{n \geq 2} B_n^{(1,0,1,0)}(x) 
\frac{ t^{n-2}}{(n-2)!} 
\end{equation}
where 
$$B_n^{(1,0,1,0)}(x) = \sum_{\sg \in S_n(1 \rightarrow n)} x^{mpp^{(1,0,1,0)}(\sg)}.$$ 
Let $S_n(i,j)$ denote the set of all permutations $\sg = \sg_1 \ldots \sg_n \in S_n$ such that $\sg _i =1$ and $\sg_j =n$ where $1 \leq i < j \leq n$.  
It is easy to see that 
for $s <i$ $\sg_s$ matches the pattern $MMP(1,0,1,0)$ in $\sg$ if and 
only if $\sg_s$ matches the pattern $MMP(0,0,1,0)$ in $\sg_1 \ldots \sg_{i-1}$.
Clearly if $i < s < j$, then $\sg_s$ matches the pattern $MMP(1,0,1,0)$ in 
$\sg$. Finally, if $j < s \leq n$, then $\sg_s$ matches the pattern $MMP(1,0,1,0)$ in $\sg$ if and 
only if $\sg_s$ matches the pattern $MMP(1,0,0,0)$ in $\sg_{j+1} \ldots 
\sg_{n}$. It follows that 
\begin{equation}\label{eq1:B1010-}
B_n^{(1,0,1,0)}(x) = \sum_{1 \leq i < j \leq n} 
\binom{n-2}{i-1,j-i-1,n-j} (j-i-1)! x^{j-i-1} R_{i-1}^{(0,0,1,0)}(x) 
R_{n-j}^{(1,0,0,0)}(x)
\end{equation}
or, equivalently, 
 \begin{eqnarray}\label{eq1:B1010}
\frac{B_n^{(1,0,1,0)}(x)}{(n-2)!} &=& \sum_{1 \leq i < j \leq n} 
x^{j-i-1} \frac{R_{i-1}^{(0,0,1,0)}(x)}{(i-1)!} 
\frac{R_{n-j}^{(1,0,0,0)}(x)}{(n-j)!} \nonumber \\
&=&\sum_{a,b,c \geq 0, a+b+c =n-2} x^{a} \frac{R_{b}^{(0,0,1,0)}(x)}{b!} 
\frac{R_{c}^{(1,0,0,0)}(x)}{c!}.
\end{eqnarray}
It is easy to see from (\ref{eq:k3}) that 
$R^{(1,0,0,0)}(t,x) = R^{(0,0,1,0)}(t,x) = (1-tx)^{-1/x}$. 
Multiplying  (\ref{eq1:B1010}) by $\frac{t^{n-2}}{(n-2)!}$ and 
summing for $n \geq 2$, we obtain 
\begin{equation}\label{eq3:B1010}
B^{(1,0,1,0)}(t,x) = \frac{1}{(1-tx)}R^{(0,0,1,0)}(t,x)R^{(0,0,1,0)}(t,x)
= (1-tx)^{-1-\frac{2}{x}}.
\end{equation}
One can use (\ref{eq3:B1010}) to show that for $n \geq 1$, 
\begin{equation}\label{eq4:B1010}
B_{2n}^{(1,0,1,0)}(x) = 2^{n-1} \prod_{i=1}^{n-1} (1+ix) 
\prod_{i=1}^{n-1} (2+(2i-1)x)
\end{equation}
and 
\begin{equation}\label{eq5:B1010}
B_{2n-1}^{(1,0,1,0)}(x) = 2^{n-1} \prod_{i=1}^{n-2} (1+ix) 
\prod_{i=1}^{n-1} (2+(2i-1)x).
\end{equation}

\begin{problem} Can one find a direct explanation of formulas (\ref{eq4:B1010}) and (\ref{eq5:B1010})? \end{problem}

For example, it follows that the number of permutations 
$\sg \in S_n(1 \rightarrow n)$ which avoid $MMP(1,0,1,0)$ is 
$2^{n-2}$. This is easy to see since to avoid 
 $MMP(1,0,1,0)$ one must have 1 and $n$ be adjacent and there 
must be decreasing sequences on either side of 1 and $n$. Thus 
such a permutation is determined by choosing for each 
$i \in \{2,\ldots, n-1\}$ whether it is placed on the left or the 
right of the adjacent pair $1n$. 

There is a second class of permutations for which 
it is easy to compute the generating function 
for the distribution of matches of the pattern 
$MMP(1,0,1,0)$. That is, suppose that $k, \ell > 0$ and $n > k+\ell$, 
Let $\beta = \beta_1 \ldots, \beta_k$ be any permutation of 
$n-k+1, \ldots, n$ and $\alpha = \alpha_1 \ldots \alpha_{\ell}$ be 
any permutation of $1, \ldots, \ell$.  Then we let 
$S_n^{\overline{\beta\alpha}}$ denote the set of permutations 
$\sg \in S_n$ such that $\beta\alpha$ is a consecutive sequence 
in $\sg$. Now suppose that $a \leq k$, $b \leq \ell$, 
$mmp^{(a,0,0,0)}(\red[\beta]) =0$,  
$mmp^{(0,0,b,0)}(\red[\alpha]) =0$, and 
\begin{equation}
R_n^{(a,0,b,0),\overline{\beta \alpha}}(x) = 
\sum_{\sg \in S_n^{\overline{\beta \alpha}}} x^{mmp^{(a,0,b,0)}(\sg)}.
\end{equation}
Then we claim that for $n > k + \ell$, 
\begin{equation}
R_n^{(a,0,b,0),\overline{\beta \alpha}}(x) = \sum_{i=0}^{n-k - \ell} 
\binom{n-k-\ell}{i} R^{(0,0,b,0)}_i(x) R^{(a,0,0,0)}_{n-k - \ell -i}(x).
\end{equation}
That is, suppose that  $\sg \in S_n^{\overline{\beta\alpha}}$ is 
a permutation where  
there are $i$ elements to the left of the occurrence of $\beta\alpha$.  
Then a $\sg_j$ with $j \leq i$ matches the pattern $MMP(a,0,b,0)$ in 
$\sg$ 
if and only if $\sg_j$ matches the pattern $MMP(0,0,b,0)$ in 
$\sg_1 \ldots \sg_i$. Similarly for $j > i+k+\ell$, 
$\sg_j$ matches the pattern $MMP(a,0,b,0)$ in $\sg$ 
if and only if $\sg_j$ matches the pattern $MMP(a,0,0,0)$ in 
$\sg_{k+\ell+1} \ldots \sg_n$. Finally our assumptions ensure 
that none of the elements 
in $\overline{\beta\alpha}$ matches the pattern 
$MMP(a,0,b,0)$ in $\sg$. It then follows 
that 
\begin{equation}
R^{(a,0,b,0),\overline{\beta\alpha}}(t,x) =
\sum_{n \geq k +\ell} \frac{t^{n-k-\ell}}{(n-k - \ell)!} 
R_n^{(a,0,b,0),\overline{\beta\alpha}}(x) = 
R^{(a,0,0,0)}(t,x) R^{(0,0,b,0)}(t,x).
\end{equation}
For example, since $R^{(0,0,1,0)}(t,x) = R^{(1,0,0,0)}(t,x) = (1-tx)^{-1/x}$, 
it follows that 
\begin{eqnarray}
R^{(1,0,1,0),\overline{n1}}(t,x) &=& \sum_{n \geq 2} \frac{t^{n-2}}{(n-2)!} 
R_n^{(1,0,1,0),\overline{n1}}(x) \nonumber \\
&=& (1-tx)^{-2/x} = 1 + \sum_{n \geq 1} \frac{t^n}{n!}  \prod_{i=0}^{n-1} 
(2+ix).
\end{eqnarray}

\section{$R_n^{(1,0,1,1)}(x)=R_n^{(0,1,1,1)}(x)=R_n^{(1,1,0,1)}(x)=R_n^{(1,1,1,0)}(x)$}

The equalities in the section title are true by Lemma \ref{sym}. We will be considering  $R_n^{(1,0,1,1)}(x)$.

In this case, we can develop a recursion 
to compute $R_n^{(1,0,1,1)}(x)$ which 
is very similar to the recursion that we developed 
to compute  $R_n^{(1,0,1,0)}(x)$.  That is, 
for any permutation $\sg = \sg_1 \ldots \sg_n \in S_n$, 
let 
\begin{eqnarray*}
C_i(\sg) &=& \chi(\sg_i \ \mbox{matches the pattern} \ MMP(1,0,1,1) \ \mbox{in} \ \sg) \ 
\mbox{and} \\
D_i(\sg) &=& \chi(\sg_i \ \mbox{matches the pattern} \ 
MMP(\emptyset,0,1,1) \ \mbox{in} \ \sg).
\end{eqnarray*}
Then let 
\begin{equation}
F^{(1,0,1,1)}_n(x_2, \ldots, x_{n-1};y_2, \ldots, y_{n-1}) = 
\sum_{\sg \in S_n} \prod_{i=2}^n x_i^{C_i(\sg)} y_i^{D_i(\sg)}.
\end{equation}
Note that $\sg_1$ and $\sg_n$  never contribute to $mmp^{(1,0,1,1)}(\sg)$ or 
$mmp^{(\emptyset,0,1,1)}(\sg)$ so that \\
$F^{(1,0,1,1)}_n(x_2, \ldots, x_{n-1};y_2, \ldots, y_{n-1})$ records all 
the information about the contributions of 
$2, \ldots, n-1$ to $mmp^{(1,0,1,1)}(\sg)$ or 
$mmp^{(\emptyset,0,1,1)}(\sg)$ as $\sg$ ranges over $S_n$.

First consider all the permutations  $\sg^{[1]}$ as $\sg$ ranges over 
$S_n$. Clearly, since each of these permutations starts with $n+1$, 
$n+1$ does not contribute to either $mmp^{(1,0,1,1)}(\sg^{[1]})$ or 
$mmp^{(\emptyset,0,1,1)}(\sg^{[1]})$. Moreover, $n+1$ does not 
effect whether any other $\sg_i$ contributes to 
$mmp^{(1,0,1,1)}(\sg^{[1]})$ or 
$mmp^{(\emptyset,0,1,1)}(\sg^{[1]})$, but it does shift the index 
over by 1. Thus it is easy to see that the contribution 
of the permutations of the form $\sg^{[1]}$ to 
$F^{(1,0,1,1)}_{n+1}(x_2,\ldots,x_{n};y_2,\ldots, y_{n})$ is just 
$F^{(1,0,1,1)}_n(x_3, \ldots, x_{n};y_3, \ldots, y_{n})$. 
Next consider all the permutations  $\sg^{[n+1]}$ as $\sg$ ranges over 
$S_n$. Clearly, since each of these permutations ends with $n+1$, 
$n+1$ does not contribute to either $mmp^{(1,0,1,1)}(\sg)$ or 
$mmp^{(\emptyset,0,1,1)}(\sg)$. However, while  $n+1$ does not 
effect whether any other $i$ contributes to $mmp^{(1,0,1,1)}(\sg)$, it 
does ensure that each $(i,\sg_i)$ that contributed to  
$mmp^{(\emptyset,0,1,1)}(\sg)$ will now contribute to 
$mmp^{(1,0,1,1)}(\sg^{[n+1]})$.
Thus it is easy to see that the contribution 
of the permutations of the form $\sg^{[n+1]}$ to 
$F^{(1,0,1,1)}_{n+1}(x_2,\ldots,x_{n};y_2,\ldots, y_{n})$ is just 
$F^{(1,0,1,1)}_n(x_2, \ldots, x_{n-1};x_2, \ldots, x_{n-1})$. 
For $2 \leq i \leq n$, again 
consider all the permutations  $\sg^{[i]}$ as $\sg$ ranges over 
$S_n$. For $j < i$, if $C_j(\sg) =1$, then 
$C_j(\sg^{[i]})  =1$ and if $D_j(\sg) =1$, then 
$C_j(\sg^{[i]})  =1$.  For $j \geq i$, if $C_j(\sg) =1$, then 
$C_{j+1}(\sg^{[i]})  =1$  and if $D_j(\sg) =1$, then 
$D_{j+1}(\sg^{[i]})  =1$. Moreover, the fact that 
$n+1$ is in position $i$ in $\sg^{[i]}$ means that 
$D_i(\sg^{[i]})  =1$. Since at most one 
of $C_i(\sg)=1$ and $D_i(\sg) =1$ holds, it follows that 
the contribution 
of the permutations of the form $\sg^{[i]}$ to 
$F^{(1,0,1,1)}_{n+1}(x_2,\ldots,x_{n};y_2,\ldots, y_{n})$ is just 
$$y_iF^{(1,0,1,1)}_n(x_2, \ldots,x_{i-1},x_{i+1},\ldots,  x_{n};x_2,\ldots,x_{i-1},y_{i+1}, \ldots, y_{n}).$$ Thus 
\begin{eqnarray}\label{1011rec}
&&F_{n+1}^{(1,0,1,1)}(x_2,\ldots, x_{n};y_2, \ldots, y_{n}) =  
F_n^{(1,0,1,1)}(x_3, \ldots, x_{n};y_3, \ldots, y_{n}) + \nonumber \\
&&F_n^{(1,0,1,1)}(x_2, \ldots, x_{n-1};x_2, \ldots, x_{n-1}) + \nonumber \\
&&  \sum_{i=2}^{n-1} y_iF_n^{(1,0,1,1)}(x_2, \ldots,x_{i-1},x_{i+1},\ldots,  x_{n};x_2,\ldots,x_{i-1},y_{i+1}, \ldots, y_{n}).
\end{eqnarray}
It is then easy to see that 
\begin{equation}\label{1011rec2}
R_n^{(1,0,1,1)}(x) = F_n^{(1,0,1,1)}(x, \ldots, x;1, \ldots,1).
\end{equation}
Using (\ref{1011rec}) and (\ref{1011rec2}) one can compute 
that \\
$R_1^{(1,0,1,1)}(x) =1$,\\
$R_2^{(1,0,1,1)}(x) =2$,\\
$R_3^{(1,0,1,1)}(x) =6$,\\
$R_4^{(1,0,1,1)}(x) =20+4x$,\\
$R_5^{(1,0,1,1)}(x) =70+42x+8x^2$,\\
$R_6^{(1,0,1,1)}(x) =252+300x+144x^2+24x^3$,\\
$R_7^{(1,0,1,1)}(x) =924+1812x+1572x^2+636x^3+96x^4$, and \\
$R_8^{(1,0,1,1)}(x) =3432+9960x+13440x^2+9576x^3+3432x^4+480x^5$.\\

In this case, the sequence 
$(R_n^{(1,0,1,1)}(0))$ is A000984 in the OEIS which 
is the sequence whose $n$-th term is the central binomial coefficient $\binom{2n}{n}$ which has lots of combinatorial interpretations. 
That is, we claim that $R_n^{(1,0,1,1)}(0)=\binom{2n-2}{n-1}$. Note that avoiding $MMP(1,0,1,1)$ is equivalent to simultaneously avoiding the patterns 1324, 1342, 2314 and 2341, which, in turn, by applying the compliment is equivalent to simultaneous avoidance of the patterns 4231, 4213, 3241 and 3214.  It is a known fact (see \cite[Subsection 6.1.1]{kit}) that the number of permutations in $S_n$ avoiding simultaneously the patterns 4132, 4123, 3124 and 3142 is given by $\binom{2n-2}{n-1}$. Mark Tiefenbruck \cite{T} has constructed 
a bijection between the set of permutations of 
$S_n$ which simultaneously avoid the patterns  4132, 4123, 3124 and 3142 
and the set of permutations of $S_n$  which simultaneously avoid the patterns 
 4231, 4213, 3241 and 3214.  Indeed,this is a special case of 
a more general bijection which will appear in a forthcoming 
paper by Remmel and Tiefenbruck \cite{RT} so we will not give 
the details here.

In this case, it is also  easy to understanding the highest 
coefficient of the polynomial $R_n^{(1,0,1,1)}(x)$. That is, 
one obtains the maximum number of occurrences of the 
pattern $MMP(1,0,1,1)$ when the permutation $\sg$ either 
starts with 1 and ends with either $2~n$ or $n~2$ or it starts 
with 2 and ends with either $1~n$ or $n~1$. Thus it is easy to 
see that the highest coefficient in $R_n^{(1,0,1,1)}(x)$ is 
$4((n-3)!)x^{n-3}$ for $n \geq 4$.

While we cannot  find the generating function 
$R^{(1,0,1,1)}(t,x)$, we can find generating functions 
for the distribution of $MMP(1,0,1,1)$-matches in 
certain restricted classes of permutations. 
That is, let 
\begin{equation}\label{B1011}
B^{(1,0,1,1)}(t,x) = \sum_{n \geq 2} B_n^{(1,0,1,1)}(x) 
\frac{ t^{n-2}}{(n-2)!} 
\end{equation}
where 
$$B_n^{(1,0,1,1)}(x) = \sum_{\sg \in S_n(1 \rightarrow n)} x^{mpp^{(1,0,1,1)}(\sg)}.$$ 
Consider the set $S_n^{(k+1)}$ of all permutations $\sg = \sg_1 \ldots \sg_n \in S_n$ such that $\sg _{k+1} =1$ and $\sg_j =n$ for some $k+1 < j \leq n$.  
It is easy to see that 
for $s \leq k$, $\sg_s$ matches the pattern $MMP(1,0,1,1)$ in $\sg$ if and 
only if $\sg_s$ matches the pattern $MMP(0,0,1,0)$ in $\sg_1 \ldots \sg_{k}$.
Clearly $\sg_{k+1} =1$ does not match the pattern $MMP(1,0,1,1)$ in $\sg$. 
If $k+1 < s \leq n$, then $\sg_s$ matches the pattern $MMP(1,0,1,1)$ in $\sg$ if and 
only if $\sg_s$ matches the pattern $MMP(1,0,0,1)$ in $\sg_{k+2} \ldots 
\sg_{n}$. It follows that 
\begin{equation}\label{eq1:B1011-}
B_n^{(1,0,1,1)}(x) = \sum_{k=0}^{n-2} 
\binom{n-2}{k} R_{k}^{(0,0,1,0)}(x) 
R_{n-k-1}^{(1,0,0,1)}(x)
\end{equation}
or, equivalently, 
 \begin{eqnarray}\label{eq1:B1011}
\frac{B_n^{(1,0,1,1)}(x)}{(n-2)!} &=& \sum_{k=0}^{n-2} 
\frac{R_{k}^{(0,0,1,0)}(x)}{k!} 
\frac{R_{n-k-1}^{(1,0,0,1)}(x)}{(n-2-k)!} 
\end{eqnarray}
Multiplying  (\ref{eq1:B1011}) by $\frac{t^{n-2}}{(n-2)!}$ and 
summing for $n \geq 2$, we obtain 
\begin{equation}\label{eq3:B1011}
B^{(1,0,1,1)}(t,x) = R^{(0,0,1,0)}(t,x)\sum_{n \geq 0} \frac{t^n}{n!} 
R_{n+1}^{(1,0,0,1)}(x).
\end{equation}
By (\ref{eq:ab00}), 
we have that 
$$\sum_{n \geq 0} \frac{t^{n}}{n!} R_{n+2}^{(1,0,0,1)}(x) = 
2 (1-tx)^{-\frac{2}{x} -1}.$$
Hence 
$$\sum_{n \geq 0} \frac{t^n}{n!} R_{n+1}^{(1,0,0,1)}(x) = 
1+\int_0^t 2 (1-zx)^{-\frac{2}{x} -1}dz  = (1-tx)^{-2/x}.$$
Since $R^{(0,0,1,0)}(t,x) = (1-tx)^{-1/x}$, 
we obtain that 
\begin{equation}\label{eq4:B1011}
B^{(1,0,1,1)}(t,x) = (1-tx)^{-3/x} = 1+3t+ \sum_{n \geq 2} 
\frac{t^n}{n!}  \prod_{i=0}^{n-1} (3+ix).
\end{equation}
For example, it follows that the number of permutations 
$\sg \in S_n(1 \rightarrow n)$ which avoid $MMP(1,0,1,1)$ is 
$3^{n}$. 

Next suppose that $k, \ell > 0$ and $n > k+\ell$, and let $\gamma_{k,\ell} = n(n-1) \ldots (n-k+1)\ell (\ell-1) \ldots 1$.  
Then we let 
$S_n^{\overline{\gamma_{k,\ell}}}$ denote the set of permutations 
$\sg \in S_n$ such that $\gamma_{k,\ell}$ is a consecutive sequence 
in $\sg$ and  
\begin{equation}
R_n^{(1,0,1,1),\overline{\gamma_{k,\ell}}}(x) = 
\sum_{\sg \in S_n^{\overline{\gamma_{k,\ell}}}} x^{mmp^{(1,0,1,1)}(\sg)}.
\end{equation}
Then it is easy to see that for $n > k + \ell$, 
\begin{equation}\label{rec1011ab}
R_n^{(1,0,1,1),\overline{\gamma_{k,\ell}}}(x) = \sum_{i=0}^{n-k - \ell} 
\binom{n-k-\ell}{i} R^{(0,0,1,0)}_i(x) R^{(1,0,0,1)}_{n-k - \ell -i}(x).
\end{equation}
That is, if $\sg \in S_n^{\overline{\gamma_{k,\ell}}}$ is such that  
there are $i$ elements to the left of the occurrence of $\gamma_{k,\ell}$, 
then a $\sg_j$ with $j \leq i$ matches the pattern $MMP(1,0,1,1)$ in 
$\sg$ 
if and only if $\sg_j$ matches the pattern $MMP((0,0,1,0)$ in 
$\sg_1 \ldots \sg_i$. Similarly for $j > i+k+\ell$, 
$\sg_j$ matches the pattern $MMP(1,0,1,1)$ in $\sg$ 
if and only if $\sg_j$ matches the pattern $MMP(1,0,0,1)$ in 
$\sg_{k+\ell+1} \ldots \sg_n$. Clearly no element that is 
part of $\gamma_{k,\ell}$ can match  $MMP(1,0,1,1)$ in $\sg$.

It then follows 
that 
\begin{equation}
R^{(1,0,1,1),\overline{\gamma_{k,\ell}}}(t,x) =
\sum_{n \geq k +\ell} \frac{t^{n-k-\ell}}{(n-k - \ell)!} 
R_n^{(1,0,1,1),\overline{\gamma_{k,\ell}}}(x) = 
R^{(1,0,0,0)}(t,x) R^{(1,0,0,1)}(t,x).
\end{equation}
Since 
$$\sum_{n \geq 1} \frac{t^{n-1}}{(n-1)!} R_n^{(1,0,0,1)}(x) = (1-tx)^{-2/x},$$
we have that 
\begin{equation}
\sum_{n \geq 0} \frac{t^n}{n!} R_n^{(1,0,0,1)}(x) = 
1+ \int_0^t (1-zx)^{-2/x}dz = \frac{1-x+(1-tx)^{1-\frac{2}{x}}}{2-x}.
\end{equation}
Hence,   
\begin{eqnarray}
R^{(1,0,1,1),\overline{n1}}(t,x) &=&
\sum_{n \geq 2} \frac{t^{n-2}}{(n-2)!} 
R_n^{(1,0,1,1),\overline{n1}}(x) \nonumber \\ 
&=& 
(1-tx)^{-1/x} \frac{1-x+(1-tx)^{1-\frac{2}{x}}}{2-x}. 
\end{eqnarray}
One can use Mathematica to show that 
\begin{eqnarray*}
&&R^{(1,0,1,1),\overline{n1}}(t,x) = 1+ 2t + (5+x)\frac{t^3}{3!} + \\
&& (14+8x +2x^2)\frac{t^4}{4!} +(41+50x+23x^2+6x^3)\frac{t^5}{5!} +\\
&&122+268x+214x^2+92x^3_24x^4)\frac{t^6}{6!} + 
(365+1283x+1689x^2+1117x^3+466x^4+120x^5)\frac{t^7}{7!} +\\
&&(1094+5660x+11412+11656x^3+6934x^4+2844x^5+720x^6)\frac{t^8}{8!} +\\
&&(3281+23524x+68042x^2+102880x^3+89849x^4+49996x^5+20268x^6+5040x^7)\frac{t^9}{9!} + \cdots
\end{eqnarray*}
The sequence $(R_n^{(1,0,1,1),\overline{n1}}(0))_{n \geq 0}$ is 
$$1,2,5,14,41,122,365,1094,3281, \ldots$$
which is sequence A007051 in the OEIS whose $n$-th term is $\frac{1+3^n}{2}$. 
This is easily explained from the recursion (\ref{rec1011ab}) 
which shows that for $n \geq 2$, 
$$R_n^{(1,0,1,1),\overline{n1}} = \sum_{k=0}^{n-2} \binom{n-2}{k} R_k^{(0,0,1,0)}(0) 
R_{n-k-1}^{(1,0,0,1)}(0).$$
But $R_k^{(0,0,1,0)} (0)=1$ since only the decreasing permutation 
$\sg = n (n-1) \ldots 2 1$ has  
no $MMP(0,0,1,0)$ matches and $R_{n-k-1}^{(1,0,0,1)}(0) =2^{n-k-1}$ 
by (\ref{eq:ab3}). 
Thus 
\begin{eqnarray*}
R_n^{(1,0,1,1),\overline{n1}}(0) &=& \sum_{k=0}^{n-2} 
\binom{n-2}{k} 2^{n-k-1}\\
&=& \frac{1}{2}( 1+ \sum_{n=0}^{n-2} \binom{n-2}{k} 2^{n-k-2}) \\
&=& \frac{1}{2}(1+(1+2)^{n-2}) = \frac{1}{2}(1+3^{n-2}).
\end{eqnarray*}
Also the coefficient of the highest term in $R_n^{(1,0,1,1),\overline{n1}}(x)$ 
is $(n-3)!$ for $n \geq 4$ which counts all the permutations that 
start with 2 and end with $n1$.

\section{The function $R_n^{(1,1,1,1)}(x)$}

We can develop a recursion 
to compute $R_n^{(1,1,1,1)}(x)$ but our recursion will 
require 4 sets of variables.   That is, 
for any permutation $\sg = \sg_1 \ldots \sg_n \in S_n$, 
let 
\begin{eqnarray*}
K_i(\sg) &=& \chi(\sg_i \ \mbox{matches the pattern} \ MMP(1,1,1,1) 
\ \mbox{in} \ \sg),\\
L_i(\sg) &=& \chi(\sg_i \ \mbox{matches the pattern} 
\ MMP(\emptyset,1,1,1) \ \mbox{in} \ \sg),\\
M_i(\sg) &=& \chi(\sg_i  \ \mbox{matches the pattern} \ 
MMP(1,\emptyset,1,1) \ \mbox{in} \ \sg), \ \mbox{and} \\
N_i(\sg) &=& \chi(\sg_i \ \mbox{matches the pattern} \ 
MMP(\emptyset,\emptyset,1,1) \ \mbox{in} \ \sg).
\end{eqnarray*}
Then let 
\begin{multline}
H^{(1,1,1,1)}_n(x_2, \ldots, x_{n-1};y_2, \ldots, y_{n-1};
z_2, \ldots, z_{n-1};w_2, \ldots, w_{n-1})\\= 
\sum_{\sg \in S_n} \prod_{i=2}^n x_i^{K_i(\sg)} y_i^{L_i(\sg)}z_i^{M_i(\sg)}
w_i^{N_i(\sg)}.
\end{multline}

First consider all the permutations  $\sg^{[1]}$ as $\sg$ ranges over 
$S_n$. The point $(1,n+1)$ in $\sg^{[1]}$  
does not contribute to either $mmp^{(1,1,1,1)}(\sg^{[1]})$, 
$mmp^{(\emptyset,1,1,1)}(\sg^{[1]})$, $mmp^{(1,\emptyset,1,1)}(\sg^{[1]})$, or 
$mmp^{(\emptyset,\emptyset,1,1)}(\sg^{[1]})$.  However,
\begin{enumerate}
\item if $K_i(\sg) =1$, then $K_{i+1}(\sg^{[1]}) =1$,
\item if $L_i(\sg) =1$, then $L_{i+1}(\sg^{[1]}) =1$,
\item if $M_i(\sg) =1$, then $K_{i+1}(\sg^{[1]}) =1$, and 
\item if $N_i(\sg) =1$, then $L_{i+1}(\sg^{[1]}) =1$.
\end{enumerate}
It follows that the contribution 
of the permutations of the form $\sg^{[1]}$ to \\
$H^{(1,1,1,1)}_{n+1}(x_2, \ldots, x_{n};y_2, \ldots, y_{n};
z_2, \ldots, z_{n};w_2, \ldots, w_{n})$ is 
$$
H^{(1,1,1,1)}_{n}(x_3, \ldots, x_{n};y_3, \ldots, y_{n};
x_3, \ldots, x_{n};y_3, \ldots, y_{n}).
$$
Next consider all the permutations  $\sg^{[n+1]}$ as $\sg$ ranges over 
$S_n$. Clearly, since each of these permutations ends with $n+1$, 
$n+1$ does not contribute to  $mmp^{(1,1,1,1)}(\sg^{[n+1]})$, 
$mmp^{(\emptyset,1,1,1)}(\sg^{[n+1]})$,  
$mmp^{(1,\emptyset,1,1)}(\sg^{[n+1]})$,or 
$mmp^{(\emptyset,\emptyset,1,1)}(\sg^{[n+1]})$.   However,
\begin{enumerate}
\item if $K_i(\sg) =1$, then $K_{i}(\sg^{[1]}) =1$,
\item if $L_i(\sg) =1$, then $K_{i}(\sg^{[1]}) =1$,
\item if $M_i(\sg) =1$, then $M_{i}(\sg^{[1]}) =1$, and 
\item if $N_i(\sg) =1$, then $M_{i}(\sg^{[1]}) =1$.
\end{enumerate}
It follows that the contribution 
of the permutations of the form $\sg^{[n+1]}$ to \\
$H^{(1,1,1,1)}_{n+1}(x_2, \ldots, x_{n};y_2, \ldots, y_{n};
z_2, \ldots, z_{n};w_2, \ldots, w_{n})$ is 
$$
H^{(1,1,1,1)}_{n}(x_2, \ldots, x_{n-1};x_2, \ldots, x_{n-1};
z_2, \ldots, z_{n-1};z_2, \ldots, z_{n-1}).
$$

For $2 \leq i \leq n$, again 
consider all the permutations  $\sg^{[i]}$ as $\sg$ ranges over 
$S_n$. First note that the fact that $n+1$ is 
in position $i$ means that $N_i(\sg^{[i]}) =1$. 
Then for $j < i$, 
\begin{enumerate}
\item if $K_j(\sg) =1$, then $K_{j}(\sg^{[i]}) =1$,
\item if $L_j(\sg) =1$, then $K_{j}(\sg^{[i]}) =1$,
\item if $M_j(\sg) =1$, then $M_{j}(\sg^{[i]}) =1$, and 
\item if $N_j(\sg) =1$, then $M_{j}(\sg^{[i]}) =1$.
\end{enumerate}
Similarly,  for $j \geq i$, 
\begin{enumerate}
\item if $K_j(\sg) =1$, then $K_{j+1}(\sg^{[i]}) =1$,
\item if $L_j(\sg) =1$, then $L_{j+1}(\sg^{[i]}) =1$,
\item if $M_j(\sg) =1$, then $K_{j+1}(\sg^{[i]}) =1$, and 
\item if $N_j(\sg) =1$, then $L_{j+1}(\sg^{[i]}) =1$.
\end{enumerate}
It follows that for $2 \leq i \leq n$, the contribution 
of the permutations of the form $\sg^{[i]}$ to \\
$H^{(1,1,1,1)}_{n+1}(x_2, \ldots, x_{n};y_2, \ldots, y_{n};
z_2, \ldots, z_{n};w_2, \ldots, w_{n})$ is 
{\tiny   
$$
w_iH^{(1,1,1,1)}_{n}(x_2, \ldots, x_{i-1},x_{i+1}, \ldots,x_{n};
x_2, \ldots, x_{i-1},y_{i+1}, \ldots,y_n;
z_2, \ldots, z_{i-1},x_{i+1}, \ldots, x_n;z_2, \ldots, z_{i-1},y_{i+1}, 
\ldots, y_n)
$$}
Thus
{\tiny \begin{eqnarray}\label{1111rec}
&&H_{n+1}^{(1,1,1,1)}(x_2,\ldots, x_{n};y_2, \ldots, y_{n};z_2, \ldots, z_{n+1};w_2, \ldots, w_{n}) =  \\
&&H^{(1,1,1,1)}_{n}(x_3, \ldots, x_{n};y_3, \ldots, y_{n};
x_3, \ldots, x_{n};y_3, \ldots, y_{n})+ \nonumber \\
&&H^{(1,1,1,1)}_{n}(x_2, \ldots, x_{n-1};x_2, \ldots, x_{n-1};
z_2, \ldots, z_{n-1};z_2, \ldots, z_{n-1})+ \nonumber \\
&&\sum_{i=2}^n w_iH^{(1,1,1,1)}_{n}
(x_2, \ldots, x_{i-1},x_{i+1}, \ldots, x_n;
x_2, \ldots, x_{i-1},y_{i+1}, \ldots,y_n;
z_2, \ldots, z_{i-1},x_{i+1}, \ldots, x_n;
z_2, \ldots, z_{i-1},y_{i+1}, \ldots, y_n). \nonumber
\end{eqnarray}}
It is then easy to see that 
\begin{equation}\label{1111rec2}
R_n^{(1,1,1,1)}(x) = F_n^{(1,1,1,1)}(x, \ldots, x;1, \ldots,1;1, \ldots,1;1, \ldots,1).
\end{equation}
Note that for all $\sg \in S_3$
$$mmp^{(1,1,1,1)}(\sg) = mmp^{(\emptyset,1,1,1)}(\sg) = 
mmp^{(1,\emptyset,1,1)}(\sg) =0.$$
However, for $mmp^{(\emptyset,\emptyset,1,1)}(\sg) =1$ 
if $\sg$ equals 132 or 231. Thus 
$$H_3^{(1,1,1,1)}(x_2,y_2,z_2,w_2)= 4+2w_2.$$ 
Using (\ref{1111rec}) and (\ref{1111rec2}), one can compute 
that \\
$R_1^{(1,1,1,1)}(x) =1$,\\
$R_2^{(1,1,1,1)}(x) =2$,\\
$R_3^{(1,1,1,1)}(x) =6$,\\
$R_4^{(1,1,1,1)}(x) =24$,\\
$R_5^{(1,1,1,1)}(x) =104+16x$,\\
$R_6^{(1,1,1,1)}(x) =464+224x+32x^2$,\\
$R_7^{(1,1,1,1)}(x) =2088+2088x+768x^2+96x^3$, and \\
$R_8^{(1,1,1,1)}(x) =9392+16096x+11056x^2+3392x^3+384x^4$.\\

In this case, the sequence $(R_n^{(1,1,1,1)}(0))_{n \geq 1}$ seems 
to be A128652 in the OEIS which is the number of square permutations 
of length $n$. There is a formula for the numbers, namely, 
$$a(n) = 2(n+2)4^{n-3}-4^{2n-5}\binom{2n-6}{n-3}.$$

\begin{problem} Can we prove this formula (directly)? \end{problem}

In this case, it is again easy to understand the highest 
coefficient of the polynomial $R_n^{(1,1,1,1)}(x)$. That is, 
one obtains the maximum number of occurrences of the 
pattern $MMP(1,1,1,1)$ when the permutation $\sg$ either \\
(i) starts with $1~(n-1)$ or $(n-1)~1$ and ends with either $2~n$ or $n~2$,\\
(ii) starts with $2~(n-1)$ or $(n-1)~2$ and ends with either $1~n$ or $n~1$,\\
(iii) starts with $1~n$ or $n~1$ and ends with either $2~(n-1)$ or $(n-1)~2$, 
or\\
(iv) starts with $2~n$ or $n~2$ and ends with either $1~(n-1)$ or $(n-1)~1$.\\
Thus it is easy to 
see that the highest coefficient in $R_n^{(1,1,1,1)}(x)$ is 
$16((n-4)!)x^{n-4}$ for $n \geq 5$.


\begin{thebibliography}{10}
\bibitem{BrCl} Petter Br\"and\'en and Anders Claesson, Mesh patterns and the expansion of permutation statistics as sums of permutation patterns, {\em Electronic J. Combin.} {\bf 18(2)} (2011), \#P5, 14pp. 
\bibitem{HilJonSigVid} \'Isak Hilmarsson, Ingibj\"org J\'onsd\'ottir, Steinunn Sigurdardottir and Sigr\'idur Vidarsd\'ottir, Wilf classification of mesh patterns of short length, in preparation.
\bibitem{kit} S. Kitaev, Patterns in permutations and words, {\em Monographs in Theoretical Computer Science} (with a foreword by Jeffrey B. Remmel), Springer-Verlag, ISBN 978-3-642-17332-5, 2011.
\bibitem{oeis} N.~J.~A.~Sloane, The on-line encyclopedia of integer sequences,
published electronically at \phantom{*} {\tt
http:/$\!\!$/www.research.att.com/\~{}njas/sequences/}.
\bibitem{T} Mark Tiefenbruck, personal communication.

\bibitem{RT} Jeffrey Remmel and Mark Tiefenbruck, Extending bijections which 
preserve marked occurrence of patterns to bijection which preserve all 
occurrences of patterns, preprint. 



\bibitem{Ulf} Henning \'Ulfarsson, A unification of permutation patterns related to Schubert varieties, arXiv:1002.4361 (2011).
\end{thebibliography}
\end{document}